\documentclass[11pt]{article}
\usepackage{mathrsfs}
\usepackage{amssymb}
\usepackage{amsfonts}
\usepackage{color,xcolor}
\usepackage{graphicx}
\usepackage{geometry}
\usepackage{manfnt}
\RequirePackage[colorlinks,citecolor=blue,urlcolor=blue]{hyperref}
\newtheorem{theorem}{Theorem}[section]
\newtheorem{corollary}{Corollary}[section]
\newtheorem{lemma}{Lemma}[section]

\newtheorem{definition}{Definition}[section]

\begin{document}
\title{\bf Well-posedness of nonlinear transport equation by stochastic perturbation }
\author{Jinlong Wei$^a$, Rongrong Tian$^b$ and Guangying Lv$^c$\thanks{Corresponding author Email: gylvmaths@henu.edu.cn}
\\ {\small \it $^a$ School of Statistics and Mathematics, Zhongnan University of} \\ {\small \it Economics and Law, Wuhan, Hubei 430073, China} \\
{\small \it $^b$ School of Mathematics and Statistics, Huazhong University of}
\\ {\small \it Science and Technology, Wuhan, Hubei 430074, China} \\ {\small \it $^c$ School of Mathematics and Statistics, Henan University}
\\ {\small \it Kaifeng, Henan 475001, China}}

\date{\today}

 \maketitle
\noindent{\hrulefill}

\vskip1mm\noindent{\bf Abstract} We are concerned with multidimensional nonlinear stochastic transport equation
driven by Brownian motions.  For irregular fluxes, by using stochastic BGK
approximations and commutator estimates, we gain the existence and uniqueness
of stochastic entropy solutions. Besides, for $BV$ initial data, the $BV_{loc}$ and
H\"{o}lder regularities are also derived for the unique stochastic entropy
solution. Particularly, we gain a regularization result, i.e. while the existence fails for the deterministic equation, we prove that a multiplicative stochastic perturbation of Brownian type is enough to render the equation well-posed. This seems to be another explicit example (the first example is given in \cite{FGP}) of a PDE of fluid dynamics that becomes well-posed under the influence of a multiplicative Brownian type noise.

 \vskip1.2mm\noindent
{\bf MSC (2010):} 60H15 (35L65 35R60)

\vskip1.2mm\noindent
{\bf Keywords:}  Transport equations; Stochastic entropy solution; Stochastic kinetic solution; Regularizing effect

 \vskip0mm\noindent{\hrulefill}
\section{Introduction}\label{sec1}\setcounter{equation}{0}
Given $T>0$, consider the following multidimensional balance law
\begin{eqnarray}\label{1.1}
\left\{
  \begin{array}{ll}
\partial_t\rho(t,x)+\mbox{div}_xF(t,x,\rho)=A(t,x,\rho), \ \ (t,x)\in (0,T)\times \mathbb{R}^d, \\
\rho(t,x)|_{t=0}=\rho_0(x), \ \  x\in   \mathbb{R}^d.
  \end{array}
\right.
\end{eqnarray}
The first complete well-posedness result for (\ref{1.1}) was settled in the pioneering paper of Kru\u{z}kov \cite{Kru}. Under the smoothness hypothesis on $F$ and $A$, Kru\u{z}kov established the existence and uniqueness of generalized solutions. For a completely satisfactory well-posedness theory of conservation laws, we refer to the monograph of Dafermos \cite{Daf}.

An important subclass of (\ref{1.1}) is that
\begin{eqnarray}\label{1.2}
\left\{
  \begin{array}{ll}
\partial_t\rho(t,x)+b(x)\cdot\nabla f(\rho)=0, \ \ (t,x)\in (0,T)\times \mathbb{R}^d, \\
\mbox{div} b=0, \ \  x\in \mathbb{R}^d, \\
\rho(t,x)|_{t=0}=\rho_0(x), \ \ x\in \mathbb{R}^d.
  \end{array}
\right.
\end{eqnarray}
It stems from the study on the complex fluid mixing in porous media flows, which occurs in a number of important scientific and technological contexts \cite{Cag,DL,GDM1,GDM2}. A particular application of (\ref{1.2}) involves two-phase fluid flow, which has been used to study the flow of water through oil in a porous medium \cite{GDM1,GDM2}.

When we deal with porous media flows, spatial variations of porous formations occur at all length scales. But only the variations at the largest length scales are reliably reconstructed from data available. The heterogeneities occurring on the smaller length scales have to be incorporated stochastically, on the basis of random fields, in geostatistical models. Consequently, the flows through such formations are stochastic \cite{GH} and in the present paper, we are interested in the stochastic perturbation given by a Brownian noise.

There are many existing papers concerning (\ref{1.2}) by a stochastic perturbation and most of them are concentrated on the following two cases:

$\bullet$ case 1: the flux is space independent, i.e. $b(x)f(\rho)$ is replaced by $F(\rho)$ (see \cite{BVW,BM,CDK,DV,FN,Hof,Kim})
\begin{eqnarray}\label{1.3}
\left\{  \begin{array}{ll}
 \partial_t\rho(t,x)+\mbox{div}_xF(\rho)
=A(t,x,\rho)\dot{W}(t),  \ (\omega,t,x)\in \Omega\times (0,T)\times \mathbb{R}^d, \\
\rho(t,x)|_{t=0}=\rho_0(x), \  x\in\mathbb{R}^d,
  \end{array}
\right.
\end{eqnarray}
where $W(t)$ is a $1$-dimensional standard Brownian motion or a cylindrical Brownian motion defined on some complete probability space ($\Omega_1, \overline{\mathcal{F}},\mathbb{P}_1,(\overline{\mathcal{F}}_{t})_{t\geq 0}$), the stochastic integration here is interpreted in It\^{o}'s.

$\bullet$ case 2: $f(\rho)=\rho$ (see \cite{AF,DR,FF,FGP,LF,Zhang})
\begin{eqnarray}\label{1.4}
\left\{
  \begin{array}{ll}
\partial_t\rho(t,x)+b(x)\cdot\nabla\rho(t,x)
+\partial_{x_i}\rho(t,x)\circ\dot{B}_i(t)=0, \ (\omega,t,x)\in \Omega\times (0,T)\times \mathbb{R}^d, \\
\rho(t,x)|_{t=0}=\rho_0(x), \  x\in\mathbb{R}^d,
  \end{array}
\right.
\end{eqnarray}
where $B(t)=(B_1(t), B_2(t), _{\cdots}, B_d(t))$ is a $d$-dimensional standard Brownian motion defined on the classical Wiener
space ($\Omega, \mathcal{F},\mathbb{P},(\mathcal{F}_{t})_{t\geq 0}$), i.e. $\Omega$ is the
space of all continuous functions from $[0,T]$ to $\mathbb{R}^d$ with
uniform convergence topology, $\mathcal{F}$ is the Borel
$\sigma$-field, $\mathbb{P}$ is the Wiener measure, $(\mathcal{F}_{t})_{t\geq 0}$ is
the natural filtration generated by the coordinate process
$B(t,\omega)=\omega(t)$. The stochastic integration with a notation
$\circ$ is interpreted in Stratonovich sense.

The various well-posedness results have been established for above two cases, we sketch some recent works, which are relevant for the present paper. Firstly, we recall some existence and uniqueness results for the weak solutions in case 1. When $W(t)$ is a $1$-dimensional Brownian motion and $A(t,x,\rho)=A(t,x)$, a change of variable reduces the equation into a hyperbolic conservation law with random flux, then using Kru\u{z}kov's method, Kim \cite{Kim} proved the existence and uniqueness of entropy solutions. For general $A$ and space-time white noise, the uniqueness theory was first founded by Feng and Nualar in \cite{FN}, but existence part was only true for $d=1$. The first complete well-posedness result for (\ref{1.3}) with cylindrical Brownian motion was obtained by Debussche and Vovelle \cite{DV} for the case of kinetic solutions. Later in \cite{Hof}, Hofmanov\'{a} extended this result and showed that the kinetic solution was the macroscopic limit of stochastic BGK approximations. Recently, by the observations that uniform spatial $BV$-bound is preserved for (\ref{1.3}) with noise functions $(A(t,x,\rho)=A(\rho))$ satisfying a Lipschitz condition, then in \cite{CDK}, Chen, Ding and Karlsen supplied a well-posedness theory of stochastic entropy solutions in $L^p$ space. Moreover, Chen, Ding and Karlsen remarked in \cite{CDK} that all the results and techniques for $A(\rho)$ can be extended easily to stochastic balance laws with additional nonhomogeneous terms, by combining with the Gr\"{o}nwall inequality, such as
\begin{eqnarray}\label{1.5}
\left\{ \begin{array}{ll}
  \partial_t\rho(t,x)+\mbox{div}_xF(x,\rho)
=A(x,\rho)\dot{W}(t)+g(x,\rho),  \ \ (\omega,t,x)\in \Omega\times (0,T)\times \mathbb{R}^d, \\
\rho(t,x)|_{t=0}=\rho_0(x), \ x\in \mathbb{R}^d,
  \end{array}
\right.
\end{eqnarray}
for a large class of non-homogeneous terms $F,g$. All above mentioned works are assumed on $\mathbb{R}^d$,  there are also many research works for bounded domain, such as \cite{BVW1,LW}.

Now, let us outline latest results for case 2.  For a bounded vector field $b$ satisfying a globally H\"{o}lder continuous and an integrability condition on the divergence, in \cite{FGP}, Flandoli, Gubinelli and Priola obtained the existence and uniqueness of $L^\infty$ solutions of (\ref{1.4}). This result was then generalized to random vector field $b$ by Duboscq and R\'{e}veillac \cite{DR}.  At the same time, in \cite{AF}, by introducing a notion of renormalized solution (see \cite{Amb,DL}), Attanasio and Flandoli extended the existence and uniqueness of weak bounded solution to the case of $BV$ vector fields. A few years later, for Lebesgue vector field $b\in L^p \ (p>d)$ and Sobolev initial data $\rho_0\in \cap _{r\geq 1}W^{1,r}$, Fedrizzi and Flandoli \cite{FF} derived the existence and uniqueness of $W^{1,r}_{loc}$ solution.

Even though the works stated above are discussed the well-posedness of stochastic weak solutions, there are some differences between in case 1 and in case 2: works stated in case 1 are focused on finding the similarity between the deterministic PDE and the stochastic PDE, and works stated in case 2 are centralized on making the difference between the deterministic PDE and the stochastic PDE.  Precisely speaking, in case 1, all the researchers paid their attention on how to extend the well-posedness theory of weak solutions from the deterministic PDE to the stochastic versions. But, in case 2, all the works are concerned with the regularization, i.e. how to make an ill-posed deterministic PDE to become well-posed by adding a noise. And, in this article, we will care for the second issue.

Compared with works between in case 1 and in case 2, we find that the multiplicative noise given in (\ref{1.3}) and (\ref{1.5}) in general does not have any regularizing effect (see \cite{CDK}), and the noise given in (\ref{1.4}) will affect the equation (see \cite{FGP}). The example given in \cite{FGP}  seems to be the first explicit example of a PDE of fluid dynamics that becomes well-posed under the influence of a (multiplicative) noise, i.e. while uniqueness may fail for the deterministic PDE, then a multiplicative stochastic perturbation of Brownian type is enough to render the equation well-posed. For the regularization, there is another very interesting result given by \cite{DFV}:
\begin{eqnarray*}
\left\{ \begin{array}{ll}
\partial_tu(t,x,v)+v\partial_xu+[E(t,x)+\epsilon \sum_{k=1}^\infty\sigma_k(x)\circ \dot{B}_{k}]\partial_vu=0, \ t>0, \ x,v\in \mathbb{R},
\\ u(t,x,v)|_{t=0}=\rho_0(x,v), \ x,v\in \mathbb{R}, \\
\rho(t,x)=\int_{\mathbb{R}}u(t,x,v)dv, \ E(t,x)=\int_{\mathbb{R}}Z(x-y)\rho(t,y)dy,
  \end{array}
\right.
\end{eqnarray*}
where $Z$ is a bounded function that is continuous everywhere except at $x=0$ and has sides limits in $0+$ and $0-$. Without the noise, the equation is well-known as the Vlasov-Poisson equation, and the solution will blow up in finite time. Under the perturbation given above, Delarue, Flandoli and Vincenzi obtained the strong unique solvability of the stochastic model for every initial configuration of distinct point charges.

Sum over the regularization of a deterministic PDE by a multiplicative noise of Brownian type, we conclude that

$\bullet$ While uniqueness fails for a deterministic PDE, then a multiplicative stochastic perturbation of Brownian type could render the equation well-posed.

$\bullet$ While global existence fails for a deterministic PDE, then a multiplicative stochastic perturbation of Brownian type could render the solution global existence.

At the same time Flandoli and his co-authors  gave some open problems. One of the open problems is how to extend the regularizing effect coming from a multiplicative stochastic perturbation of Brownian type to the nonlinear conservation law (\ref{1.1}) ? Currently, we concrete this open problem as the following:

$\bullet$ Does there exist a deterministic PDE, such that while existence fails, then a multiplicative stochastic perturbation of Brownian type renders the equation well-posed ?

Inspired by the papers \cite{CDK,DFV,FGP}, in this paper we consider the following Cauchy problem
\begin{eqnarray}\label{1.6}
\left\{
  \begin{array}{ll}
\partial_t\rho(t,x)+b(x)\cdot\nabla f(\rho)
+\partial_{x_i}\rho(t,x)\circ\dot{B}_i(t)=0, \ \ (\omega,t,x)\in \Omega\times (0,T)\times \mathbb{R}^d, \\
\rho(t,x)|_{t=0}=\rho_0(x), \ x\in \mathbb{R}^d,
  \end{array}
\right.
\end{eqnarray}
here the choice of Stratonovich integral in (\ref{1.6}) is motivated by hyperbolicity as well. Indeed, if $b$ and $f$ are regular, $\rho_0$ smooth, It\^{o}'s formula implies that $u(t,x)=\rho(t,x+B(t))$ satisfies a nonlinear transport equation with random vector field:
\begin{eqnarray*}
\left\{
  \begin{array}{ll}
\partial_tu(t,x)+b(x+B(t))\cdot\nabla f(u)=0, \ \ (t,x)\in (0,T)\times \mathbb{R}^d, \\
u(t,x)|_{t=0}=\rho_0(x), \  x\in \mathbb{R}^d.
  \end{array}
\right.
\end{eqnarray*}
Hence, if the deterministic equation is hyperbolic, so does the stochastic equation. In this sense, the Stratonovich integral keep the original structure.

Under some suitable assumptions, we will prove the existence and uniqueness of stochastic entropy solutions for the Cauchy problem (\ref{1.6}). Moreover, for $BV$ initial data, some regularity results for (\ref{1.6}) are established.
The well-posedness as well as the regularity of solutions for the general case ($b(x)\cdot\nabla f(\rho)$ is replaced by $b(x,\rho)\cdot\rho$ in (\ref{1.6})) will be discussed in our further paper. In particular, these results (existence, uniqueness and regularity) are applicable to the stochastic transport equation (\ref{1.4}). However, for the deterministic transport equation, we can construct a counterexample for $d\geq 2$ on the existence for such
$BV_{loc}$-solutions. In this sense, the noise has a regularizing effect and we give a positive answer for above open problem. This seems to be an explicit example of a PDE of fluid dynamics that the existence of the weak solution fails, while a stochastic perturbation of Brownian type is enough to render the equation well-posed.

Now, we explain how we can solve the problem. In papers \cite{DFV} and \cite{FGP}, the authors used the idea of stochastic flow coupled with a PDE. The method is available to the linear case, but not available to the nonlinear case. That is, it is hard to generalize the results to the nonlinear case by using the method of \cite{DFV} and \cite{FGP}. Fortunately, Lions, Perthame and Tadmor \cite{LPT} introduced the kinetic formulation method and the advantage of the method is the nonlinear PDE can be changed into
a linear PDE. Thus we can combine the two methods together to solve the problem. Comparing (\ref{1.6}) with (\ref{2.2}), it is easy to find that one can change the nonlinear term $b(x)\cdot\nabla f(\rho)$ in (\ref{1.6}) into a linear term $f'(v)b(x)\cdot\nabla u$. In this paper, we only consider the special case (\ref{1.6}) because this equation has background in physics. Meanwhile, we remark that our proofs are different from that in \cite{DFV} and \cite{FGP}.
We also note that our results and the method used here are different from those of \cite{FG}. In \cite{FG}, the author assumed that the initial data are sufficient regularity, see
\cite[Theorem 1.1]{FG}.

The rest of the paper is structured as follows. In Section \ref{sec2}, we introduce some notions for solutions to (\ref{1.6}),
give some preliminaries and state out the main results: Theorems \ref{the2.1}
and \ref{the2.2}. In Sections \ref{sec3}-\ref{sec4}, the proof of
Theorem \ref{the2.1} is given. Section \ref{sec3} is concerned with the existence of the stochastic entropy solution and the proof of the uniqueness is given in Section \ref{sec4}. Section 5 is devoted to prove Theorem \ref{the2.2}. In last section, we construct a counterexample to illustrate the non-existence on $BV_{loc}$ solutions first, and then give some concluding remarks.

We end up the section by introducing some notations.

\textbf{Notations.} ${\mathcal D}(\mathbb{R}^d)$ and ${\mathcal D}(\mathbb{R}^{d+1})$ stand for the set of all smooth functions on $\mathbb{R}^d$ and $\mathbb{R}^{d+1}$ with compact supports respectively. Correspondingly, ${\mathcal D}_+(\mathbb{R}^d)$ and ${\mathcal D}_+(\mathbb{R}^{d+1})$ represents the non-negative elements in ${\mathcal D}(\mathbb{R}^d)$ and ${\mathcal D}(\mathbb{R}^{d+1})$ respectively. For an open subset $D$ in $\mathbb{R}^{d+1}$, ${\mathcal D}(D)$ and ${\mathcal D}_+(D)$ can be defined similarly. $\langle \ , \ \rangle_v$ denotes the duality between ${\mathcal D}(\mathbb{R})$ and ${\mathcal D}^\prime(\mathbb{R})$. Let $T>0$ be a given real number, $\langle \ , \ \rangle_{t,x,v}$  is the duality between ${\mathcal D}([0,T)\times\mathbb{R}^{d+1})$ and ${\mathcal D}^\prime([0,T)\times\mathbb{R}^{d+1})$. $x\in \mathbb{R}^d$ is always assumed. $\mathbb{R}^{d+1}_{x,v}=\{(x,v),  x\in\mathbb{R}^d, v\in\mathbb{R}\}$, $L^\infty_{\omega,t,x}:=L^\infty(\Omega\times [0,T]\times\mathbb{R}^d)$,
${\mathcal C}_t(L^1_{\omega,x}):={\mathcal C}([0,T];L^1(\Omega\times\mathbb{R}^d))$, $BV_x=BV(\mathbb{R}^d)$ and other notations are similar. Here a function $\zeta\in BV_x$ means $\zeta\in L^1_x$ and $D\zeta$ is a finite measure on $\mathbb{R}^d$. ${\mathcal M}_b([0,T]\times\mathbb{R}^{d+1}_{x,v})$ is the space of bounded non-negative measures over $[0,T]\times\mathbb{R}^{d+1}_{x,v}$, with norm given by the total variation of measures, corresponding  ${\mathcal M}_b(K)$ is the space of bounded non-negative measures over $K\subset\subset\mathbb{R}^{d+1}_{x,v}$. Given a measurable function $\varsigma$, $\varsigma^+$ is defined by $\max\{\varsigma,0\}$ and $\varsigma^-$ is $\max\{-\varsigma,0\}$. For $r\in \mathbb{R}$, $\mbox{sgn}(r)=1_{(0,\infty)}(r)-1_{(-\infty,0)}(r)$. For every $R>0$, $B_R:=\{x\in{\mathbb R}^d:|x|<R\}$. Almost surely can be abbreviated to $a.s.$. ${\mathbb N}$ is natural numbers set. The summation convention is enforced throughout this article.

\section{Preliminaries and main results}\label{sec2}
\setcounter{equation}{0}
\subsection{Definitions}\label{sec2.1}

\begin{definition} \label{def2.1} Let ($\Omega, \mathcal{F},{\mathbb P},(\mathcal{F}_{t})_{t\geq 0}$) be a stochastic basis, described in the introduction. Given $T>0$ and a random process (or random field) $\{\zeta_t\}_{0\leq t\leq T}$. Let ${\mathcal P}$ denote the smallest $\sigma$-field of subsets of $[0,T]\times\Omega$, which is generated by subsets of $[0,T]\times\Omega$ of the form $(s,t]\times E_1$ with $E_1 \in \mathcal{F}_s$ for $0 \leq s<t\leq T$ and $\{0\}\times E_2$ with $E_2\in \mathcal{F}$. The random process (or random field) $\{\zeta_t\}_{0\leq t\leq T}$ is said to be predictable if the function $(t,\omega)\rightarrow \zeta_t(\omega)$ is ${\mathcal P}$-measurable on $[0,T]\times\Omega$, and the random process (or random field) $\{\zeta_t\}_{0\leq t\leq T}$ is said to be $\{\mathcal{F}_{t}\}_{0\leq t\leq T}$-adapted if for every $t\in [0,T]$, the random variable (or random field) $\zeta_t$ is $\mathcal{F}_{t}$-measurable.
\end{definition}

Given $s\in[0,T]$ and $x\in \mathbb{R}^d$, consider the stochastic differential equation (SDE for short) in
$\mathbb{R}^d$:
\begin{eqnarray}\label{2.1}
dX(s,t)=a(X(s,t))dt+dB(t), \ \ t\in[s,T ], \ \ X(s,t)|_{t=s}=x.
\end{eqnarray}

\begin{definition} [\cite{Kun}, P114] \label{def2.2} A stochastic homeomorphism flow (of class  ${\mathcal C}^{1,\alpha}$) on the stochastic basis $(\Omega, \mathcal{F},{\mathbb P}, (\mathcal{F}_t)_{0\leq t\leq T})$ associated to (\ref{2.1}) is a
map $(s,t,x,\omega) \rightarrow X(s,t,x)(\omega)$, defined for
$0\leq s \leq t \leq T, \ x\in \mathbb{R}^d, \ \omega \in \Omega$ with
values in $\mathbb{R}^d$, such that

(i) given any $s\in [0,T],\  x \in \mathbb{R}^d$, the process $\{X(s,\cdot,x)\}= \{X(s,t,x), \ t\in [s,T]\}$ is a continuous $\{\mathcal{F}_{s,t}\}_{s\leq t\leq T}$-adapted solution of (\ref{2.1});

(ii) ${\mathbb P}-a.s.$, for all  $0\leq s \leq t \leq T$, the functions $X(s,t,x), \ X^{-1}(s,t,x)$ are continuous in $(s,t,x)$;

(iii) ${\mathbb P}-a.s., \   X(s,t,x)=X(r,t,X(s,r,x))$  for all
$0\leq s\leq r \leq t \leq T$, $x\in \mathbb{R}^d$ and $X(s,s,x)=x$.
\end{definition}

Now, we give the notion of stochastic entropy solution for (\ref{1.6}).

\begin{definition} \label{def2.3} Let $\rho$ be a predictable random field, which lies in $L^\infty_{\omega,t,x}\cap {\mathcal C}_t(L^1_{\omega,x})$. We call $\rho$ a stochastic weak solution of (\ref{1.6}) if for every $\varphi\in {\mathcal D}(\mathbb{R}^d)$ and every $t\in [0,T)$,
\begin{eqnarray*}
\int\limits_{\mathbb{R}^d}\varphi(x)\rho(t,x)dx&=& \int\limits_{\mathbb{R}^d}\varphi(x)\rho_0(x)dx+
\int\limits^t_0\int\limits_{\mathbb{R}^d}f(\rho)\mbox{div}(b(x)\varphi(x))dxds
\cr\cr&&+\int\limits^t_0\circ
dB_i(s)\int\limits_{\mathbb{R}^d}\partial_{x_i}\varphi(x)\rho(s,x)dx,  \quad
{\mathbb P}-a.s..
\end{eqnarray*}
The stochastic weak solution is said to be a stochastic entropy solution, if for every $\eta\in\Xi$,
\begin{eqnarray*}
\partial_t\eta(\rho)+b(x)\cdot\nabla Q(\rho)+
\partial_{x_i}\eta(\rho) \circ \dot{B}_i(t) \leq 0, \ \
{\mathbb P}-a.s.,
\end{eqnarray*}
in the sense of distributions, i.e. for every $\varphi\in {\mathcal D}_+(\mathbb{R}^d)$, $\psi\in {\mathcal D}_+([0,T))$
\begin{eqnarray*}
&&\int\limits_0^T\int\limits_{\mathbb{R}^d}\Psi(t,x)\eta(\rho(t,x))dxdt+
\int\limits^T_0\int\limits_{\mathbb{R}^d}Q(\rho)\mbox{div}(b(x)\Psi(t,x))dxdt
\cr\cr&&+ \int\limits_{\mathbb{R}^d}\Psi(0,x)\eta(\rho_0(x))dx+\int\limits^T_0\circ
dB_i(t)\int\limits_{\mathbb{R}^d}\partial_{x_i}\Psi(t,x)\eta(\rho(t,x))dx\geq 0,  \
{\mathbb P}-a.s.,
\end{eqnarray*}
where $\Phi=\varphi\psi$ and
$$
Q(\rho)=\int\limits^\rho \eta^\prime(v)f^\prime(v)dv, \ \Xi=\{c_0\rho+\sum_{k=1}^nc_k|\rho-\rho_k|, \ c_0, \rho_k, c_k\in
\mathbb{R} \ are \ constants, \ n\in {\mathbb N} \}.
$$
\end{definition}

Our proofs for the existence and uniqueness of stochastic entropy solutions is based upon rewriting (\ref{1.6}) in its kinetic form using the classical Maxwellian
\begin{eqnarray}\label{2.2}
\left\{\begin{array}{ll}
\partial_tu(t,x,v)+f^\prime(v)b(x)\cdot\nabla_xu
+\partial_{x_i}u\circ \dot{B}_i=\partial_vm, \ \ (\omega,t,x,v)\in
\Omega\times(0,T)\times \mathbb{R}^{d+1}_{x,v}, \\
u(t,x,v)|_{t=0}=\chi_{\rho_0(x)}(v), \  (x,v)\in \mathbb{R}^{d+1}_{x,v},
\end{array}
\right.
\end{eqnarray}
where $u=u(t,x,v)=\chi_{\rho(t,x)}(v)=1_{(0,\rho(t,x))}(v)-1_{(\rho(t,x),0)}(v)$, $0\leq m$ is a bounded and predictable measure, supported in $\Omega\times(0,T)\times \mathbb{R}^d\times[-N,N]$ with $N=\|\rho\|_{L^\infty_{\omega,t,x}}$. Moreover $m$ has a continuous version. Here the boundedness is understood that ${\mathbb E} m([0,T)\times\mathbb{R}^{d+1}_{x,v})<\infty$, the predictability and continuity of the measure $m$ are interpreted as follows: for every $\varphi\in {\mathcal D}(\mathbb{R}^d), \psi\in {\mathcal D}(\mathbb{R})$,
\begin{eqnarray}\label{2.3}
{\Big \{ } \int\limits_0^t\int\limits_{\mathbb{R}^{d+1}_{x,v}}\varphi(x)\psi(v)m(dx,dv,dr){\Big \}}_{0\leq t< T} \ \mbox{is predictable}
\end{eqnarray}
and
\begin{eqnarray}\label{2.4}
\lim_{s\rightarrow t}\int\limits_0^s\int\limits_{\mathbb{R}^{d+1}_{x,v}}\varphi(x)\psi(v)m(dx,dv,dr)=
\int\limits_0^t\int\limits_{\mathbb{R}^{d+1}_{x,v}}\varphi(x)\psi(v)m(dx,dv,dr), \ \mbox{for every} \ 0\leq t<T.
\end{eqnarray}

Then, with the help of  the equivalence between the stochastic entropy solution of (\ref{1.6}) and the stochastic weak solution of (\ref{2.2}) (see Lemma \ref{lem2.1}), we establish the well-posedness of (\ref{1.6}). Here the notion of stochastic weak solution for (\ref{2.2}) is given below:

\begin{definition}\label{def2.4} Assume that $0\leq m(t,x,v)$ is a bounded, predictable and continuous (in time) measure, supported in $v$. Let $u$ be a predictable random field, which lies in $L^\infty_{\omega,t,x}\cap {\mathcal C}_t(L^1_{\omega,x,v})\cap L^\infty_{\omega,t,x}(L^1_v)$. We call $u$ a stochastic weak solution of (\ref{2.2}) if
for every $\varphi\in {\mathcal D}(\mathbb{R}^d)$, every $\psi\in {\mathcal D}(\mathbb{R})$, and every $t\in [0,T)$,
\begin{eqnarray}\label{2.5}
&&\int\limits_{\mathbb{R}^{d+1}_{x,v}}\varphi(x)\psi(v)u(t,x,v)dxdv\cr\cr&=& \int\limits_{\mathbb{R}^{d+1}_{x,v}}\varphi(x)\psi(v)\chi_{\rho_0(x)}(v)dxdv+
\int\limits^t_0\int\limits_{\mathbb{R}^{d+1}_{x,v}}f^\prime(v)\mbox{div}(b(x)\varphi(x))\psi(v)
u(s,x,v)dxdvds\cr\cr&&+\int\limits^t_0\circ
dB_i(s)\int\limits_{\mathbb{R}^{d+1}_{x,v}}\psi\partial_{x_i}\varphi udxdv
-\int\limits^t_0\int\limits_{\mathbb{R}^{d+1}_{x,v}}\partial_v\psi\varphi m(ds,dx,dv),  \quad
{\mathbb P}-a.s..
\end{eqnarray}
\end{definition}

\subsection{Useful Lemmas}\label{sec2.2}
 Initially, let us establish the equivalence between the
 stochastic entropy solution to (\ref{1.6}) and the stochastic weak solution to (\ref{2.2}).
\begin{lemma} \label{lem2.1} \textbf{(Stochastic kinetic formulation)} Let $\rho_0\in L^\infty\cap L^1(\mathbb{R}^d)$, $b\in {\mathcal C}_b(\mathbb{R}^d;\mathbb{R}^d), \ \mbox{div} b\in L^\infty(\mathbb{R}^d),\ f\in {\mathcal C}^1(\mathbb{R})$ and $u=\chi_{\rho}(v)$. $0\leq m(t,x,v)$ is a bounded, predictable and continuous (in time) measure, supported in $v$.

(i) If $\rho$ is a stochastic entropy solution of (\ref{1.6}), then $u$ is a stochastic weak solution of (\ref{2.2}).

(ii) Conversely, if $u$ is a stochastic weak solution of (\ref{2.2}), then $\rho$ is a stochastic entropy solution of (\ref{1.6}).
 \end{lemma}
\textbf{Proof.} For every $\alpha_1, \alpha_2\in\mathbb{R}$,
$$
\int\limits_{\mathbb{R}}|\chi_{\alpha_1}(v)-\chi_{\alpha_2}(v)|dv=|\alpha_1-\alpha_2|,
$$
so $\rho\in L^\infty_{\omega,t,x}\cap {\mathcal C}_t(L^1_{\omega,x})$  implies $u\in {\mathcal C}_t(L^1_{\omega,x,v})\cap L^\infty_{\omega,t,x}(L^1_v)$, and vice versa. We need to survey the rest of Lemma \ref{lem2.1}.

(i) Assume that $\rho$ is a stochastic entropy solution of
(\ref{1.6}). Then, for every $v\in \mathbb{R}$, it renders that
\begin{eqnarray}\label{2.6}
\partial_t\eta(\rho,v)+b(x)\cdot\nabla_xQ(\rho,v)+
\partial_{x_i}\eta(\rho,v) \circ \dot{B}_i(t)=-2m, \quad {\mathbb P}-a.s.,
\end{eqnarray}
where
\begin{eqnarray}\label{2.7}
\left\{
  \begin{array}{ll}
    \eta(\rho,v)=|\rho-v|-|v|, \\
Q(\rho,v)=\mbox{sgn}(\rho-v)[f(\rho)-f(v)]-\mbox{sgn}(v)f(v), \\
 m \ \mbox{is a nonnegative measure on} \
\Omega\times[0,T]\times \mathbb{R}^{d+1}_{x,v}.
  \end{array}
\right.
\end{eqnarray}

If one differentiates $\eta$, $Q$ and $m$ in $v$ in distributional sense, then
\begin{eqnarray}\label{2.8}
\partial_v\eta(\rho,v)=-2u(t,x,v),  \ \
\partial_vQ(\rho,v)=-2f^\prime(v)u(t,x,v), \quad {\mathbb P}-a.s..
\end{eqnarray}
Thus one derives the identity (\ref{2.2}) in the sense of distributions for almost all $\omega\in\Omega$.

In view of that $\rho$ is bounded, from (\ref{2.6})
and (\ref{2.7}), $m$ is supported in $[0,T]\times \mathbb{R}^d\times [-N,N]$, with $N=\|\rho\|_{L^\infty_{\omega,t,x}}$. The predictability can be derived from the identity (\ref{2.5}), accordingly,
it remains to examine that $m$ is bounded and continuous in $t$.

Since $m\geq 0$ and it is supported in a compact subset in $v$, we
obtain
\begin{eqnarray}\label{2.9}
0&\leq&\langle m, \ \Psi\otimes 1 \rangle_{t,x,v} \cr\cr&=&
-\langle\partial_tu+f^\prime(v)b(x)\cdot \nabla_{x}
u+\partial_{x_i}u\circ \dot{B}_i(t), \ \Psi\otimes v
\rangle_{t,x,v}\cr\cr
&=&-\langle\partial_tu+f^\prime(v)b(x)\cdot \nabla_{x}
u+\partial_{x_i}u\dot{B}_i(t)-\frac{1}{2}\Delta_x u(t,x,v), \ \Psi\otimes v \rangle_{t,x,v},
\end{eqnarray}
for every $\Psi\in {\mathcal D}_+([0,T)\times \mathbb{R}^d)$, and for almost all
$\omega\in\Omega$.

Thanks to the fact
\begin{eqnarray*}
\int\limits_{\mathbb{R}}g^\prime(v)u(t,x,v)dv=g(\rho(t,x))-g(0), \ \mbox{for
every} \ g\in W^{1,1}_{loc}({\mathbb{R}}),
\end{eqnarray*}
one computes from (\ref{2.9}) that
\begin{eqnarray}\label{2.10}
&& \langle m, \ \Psi\otimes 1 \rangle_{t,x,v}\cr\cr&=&\frac{1}{2}\int\limits^T_0\int\limits_{
\mathbb{R}^d}\partial_t\Psi(t,x)\rho^2(t,x)dxdt+\frac{1}{2}\int\limits_{\mathbb{R}^d}
\Psi(0,x)\rho^2_0(x)dx\cr \cr &&
+\int\limits^T_0\int\limits_{\mathbb{R}^d}\Big[\rho(t,x)
f(\rho(t,x))-\int\limits^{\rho(t,x)}_0f(v)dv\Big] \mbox{div}_x(b(x)\Psi(t,x))
dxdt \cr \cr && + \frac{1}{4}\int\limits^T_0\int\limits_{\mathbb{R}^d} \rho^2(t,x)\partial^2_{x_i,x_j}\Psi(t,x)dxdt
+ \frac{1}{2}\int\limits^T_0\int\limits_{\mathbb{R}^d}\rho^2(t,x)
\partial_{x_i}\Psi(t,x)dxdB_i(t)
 \cr \cr &\leq&
\frac{1}{2}\int\limits^T_0\int\limits_{\mathbb{R}^d}
\partial_t\Psi\rho^2(t,x)dxdt+\frac{1}{2}
\int\limits_{\mathbb{R}^d}\Psi(0,x)\rho^2_0(x)dx+
C\int\limits^T_0\int\limits_{\mathbb{R}^d}|\rho(t,x)|\Big[\Psi(t,x)+|\nabla_x
\Psi(t,x)|\Big] dxdt \cr \cr
&& + \frac{1}{4}\int\limits^T_0\int\limits_{\mathbb{R}^d}\rho^2(t,x)
\partial^2_{x_i,x_j} \Psi(t,x)|dxdt + \frac{1}{2}\int\limits^T_0\int\limits_{\mathbb{R}^d}\rho^2(t,x)\partial_{x_i}
\Psi(t,x)dxdB_i(t), \quad {\mathbb P}-a.s..
\end{eqnarray}

Obviously, (\ref{2.10}) holds ad hoc for
$\psi(t,x)=\psi(t)\theta_n(x)$, where $\psi\in {\mathcal D}_+([0,T)), \
\theta\in {\mathcal D}_+(\mathbb{R}^d)$,
\begin{eqnarray}\label{2.11}
\theta_n(x)=\theta(\frac{x}{n}),
 \ \theta(x)=\cases{1, \ \mbox{when} \
|x|\leq 1, \cr 0, \ \mbox{when} \ |x|> 2.}
\end{eqnarray}
For this fixed $n$, by an approximation demonstration, one can
fetch
\begin{eqnarray}\label{2.12}
\psi(t)=\cases{ \ \ \ \ \ \ \ 1, \ \ \ \ \ \ \
t\in[0,T_0-\frac{1}{n}], \cr -n(t-T_0), \ t\in (T_0-\frac{1}{n},T_0], \cr
\ \ \ \ \ \ \ 0,   \ \ \ \ \ \ \   t\geq T_0,}
\end{eqnarray}
with every $T_0\in [0,T)$.

Observing that the It\^{o} isometry,
\begin{eqnarray}\label{2.13}
{\mathbb E}{\Big[} \int\limits^{T_0}_0\int\limits_{\mathbb{R}^d}\rho^2(t,x)\partial_{x_i}(\theta_n(x))\psi(t)dxdB_i(t){\Big]}^2=
\frac{1}{n^2}{\mathbb E} \int\limits^{T_0}_0 {\Big[}
\int\limits_{\mathbb{R}^d}\rho^2(t,x)|\partial_{x_i}\theta_n(x)|\psi(t)dx{\Big]}^2dt,
\end{eqnarray}

Thus we gain from (\ref{2.10})-(\ref{2.13}) by letting
$n\rightarrow\infty$,
\begin{eqnarray}\label{2.14}
&&\int\limits_0^{T_0}\int\limits_{\mathbb{R}^{d+1}_{x,v}}m(dt,dx,dv) \cr\cr
&\leq&\frac{1}{2}{\Big[}\int\limits_{\mathbb{R}^d}\rho^2_0(x)dx -\int\limits_{
\mathbb{R}^d}\rho^2({T_0},x)dx {\Big]}+ C\int\limits^{T_0}_0\int\limits_{
\mathbb{R}^d}|\rho(t,x)|dxdt
\cr\cr
&\leq&\frac{1}{2}\int\limits_{\mathbb{R}^d}\rho^2_0(x)dx + C\int\limits^T_0\int\limits_{\mathbb{R}^d}|\rho(t,x)|dxdt,\quad {\mathbb P}-a.s..
\end{eqnarray}
Hence $m$ is bounded on $\Omega\times[0,T)\times \mathbb{R}^{d+1}_{x,v}$ and $m([0,T)\times \mathbb{R}^{d+1}_{x,v})\in L^1(\Omega)$.

Specially, when $T_0\rightarrow 0$, from the first inequality in (\ref{2.14}) we obtain
$$
\lim_{T_0\rightarrow 0}\int\limits_0^{T_0}\int\limits_{\mathbb{R}^d}\int\limits_{
\mathbb{R}}m(dt,dx,dv)=0, \quad {\mathbb P}-a.s..
$$

The arguments employed above for $0$ and $T_0$ adapted to any $0\leq
s, t<T$ now, yields that
$$
\lim_{t\rightarrow s}\int\limits_s^t\int\limits_{\mathbb{R}^{d+1}_{x,v}}m(dr,dx,dv)=
0,
$$
which hints $m$ is continuous in $t$. Therefore $u$ is a
stochastic weak solution of (\ref{2.2}).

(ii) Let us show the reverse fact. Given $\epsilon>0$ and
$\bar{\rho}\in \mathbb{R}$, set
$$
\eta_\epsilon(\cdot, \bar{\rho})=(\sqrt{(\cdot-\bar{\rho})^2+
\epsilon^2}-\epsilon)-|\bar{\rho}|\in {\mathcal C}^2(\mathbb{R}),
$$
then $\eta_\epsilon$ is convex, $\eta_\epsilon^\prime(\cdot,
\bar{\rho})\in {\mathcal C}_b(\mathbb{R})$, and
$$
\eta_\epsilon(r, \bar{\rho})\longrightarrow |r-\bar{\rho}|
-|\bar{\rho}| \ \ \mbox{as} \ \ \epsilon \longrightarrow 0.
$$

In a consequence of $u(t,x,v)$ solving (\ref{2.2}), it follows that
\begin{eqnarray}\label{2.15}
&&\langle \partial_vm, \ \psi\eta_\epsilon^\prime(v,
\bar{\rho})\xi_{k}(v) \rangle_{t,x,v}\cr\cr&=& \langle
\partial_t u+f^\prime(v)b(x)\cdot \nabla_{x} u+
\partial_{x_i}u\circ \dot{B}_i(t) , \ \psi\eta_\epsilon^\prime(v,
\bar{\rho})\xi_{k}(v) \rangle_{t,x,v},
\end{eqnarray}
for every $\Psi\in {\mathcal D}_+([0,T)\times \mathbb{R}^d)$, $\xi\in
{\mathcal D}_+(\mathbb{R})$, where
\begin{eqnarray}\label{2.16}
\xi_{k}(v)=\xi(\frac{v}{k}), \  0\leq \xi\leq 1,  \ \xi(v)=\cases{1,
\ \ \mbox{when} \ \ |v|\leq 1, \cr 0,  \ \ \mbox{when} \ \ |v|\geq
2.}
\end{eqnarray}
Applying the partial integration, and noticing that $m$ is nonnegative, one deduces
\begin{eqnarray}\label{2.17}
\lim_{k\rightarrow \infty}\langle \partial_vm, \
\Psi\eta_\epsilon^\prime(v, \bar{\rho})\xi_{k} \rangle_{t,x,v}
=-\lim_{k\rightarrow \infty}\langle m, \
\Psi[\eta_\epsilon^{\prime\prime}(v,
\bar{\rho})\xi_{k}+\eta_\epsilon^\prime(v,
\bar{\rho})\xi^\prime_{k}] \rangle_{t,x,v}\leq 0,  \ \ {\mathbb P}-a.s..
\end{eqnarray}

Upon using (\ref{2.17}), from
(\ref{2.15}), we derive
\begin{eqnarray}\label{2.18}
&&\int\limits^T_0dt\int\limits_{\mathbb{R}^d}\partial_t\Psi(t,x)[\eta_\epsilon(\rho,\bar{\rho}
)-\eta_\epsilon(0,\bar{\rho})]dx
+\int\limits^T_0dt\int\limits_{\mathbb{R}^d}Q_\epsilon(\rho,\bar{\rho})\mbox{div}_x(b(x)\Psi)
dx \cr \cr &\geq&-
\int\limits_{\mathbb{R}^d}\Psi(0,x)[\eta_\epsilon(\rho_0,\bar{\rho})-\eta_\epsilon(0,\bar{\rho})]dx
- \int\limits^T_0\circ dB_i\int\limits_{\mathbb{R}^d}\partial_{x_i}\Psi
\eta_\epsilon(\rho,\bar{\rho})dx,
\end{eqnarray}
by taking $k$ to infinity, here
\begin{eqnarray*}
 Q_\epsilon(\rho,\bar{\rho})=\int\limits_{\mathbb{R}}f^\prime(v)
\eta_\epsilon^\prime(v,\bar{\rho})u(t,x,v)dv.
\end{eqnarray*}

On the other hand
$$
\lim_{\epsilon\rightarrow
0}\eta_\epsilon^\prime(v,\bar{\rho})=\mbox{sign}(v-\bar{\rho})
$$
and
$$
\lim_{\epsilon\rightarrow
0}Q_\epsilon(\rho,\bar{\rho})=\mbox{sign}(\rho-\bar{\rho})
[f(\rho)-f(\bar{\rho})]-\mbox{sign}\bar{\rho}[f(\bar{\rho})-f(0)],
$$
for a.s. $(\omega,t,x)\in\Omega\times [0,T)\times \mathbb{R}^d$.

If one lets $\epsilon$ approach to zero in (\ref{2.18}), we
attain the inequality in Definition \ref{def2.3}, thus $\rho$ is a stochastic entropy solution. $\Box$

\vskip2mm\noindent
\textbf{Remark 2.1.} Without the noise and the flux is assumed to be $x$-independent, the equivalence between (\ref{1.6}) and (\ref{2.2}) is the classical kinetic formula (see Theorem 1 in \cite{LPT}),  we extend the result in \cite{LPT} to the case of (\ref{1.6}) and (\ref{2.2}) here. Some other various extensions to hyperbolic-parabolic equations one consults to \cite{CP,KU}.

For future use, we review some other lemmas.
\begin{lemma} [\cite{RY} P29] \label{lem2.2} \textbf{(L\'{e}vy's modulus of continuity)} Let $B(t)$
be the $d$-dimensional Brownian motion (described in introduction), then
\begin{eqnarray}\label{2.19}
{\mathbb P}\Big[\limsup_{\epsilon\downarrow0}\sup_{0\leq t_1<t_2\leq T,t_2-t_1<\epsilon}\frac{|B(t_2)-B(t_1)|}{\sqrt{2\epsilon\log\log(1/\epsilon)}}=\sqrt{d}\Big]=1.
\end{eqnarray}
 \end{lemma}

\begin{lemma} [\cite{Amb}] \label{lem2.3}  \textbf{(i) (Commutator Estimate)}
 Assume $b\in BV_{loc}(\mathbb{R}^d;\mathbb{R}^d)$, $\mbox{div} b\in L^\infty(\mathbb{R}^d)$. For every $w\in L^\infty_{loc}(\mathbb{R}^d)$ and for every smooth and even convolution kernel $\varrho_1: \mathbb{R}^d
\mapsto[0,\infty)$ with $\mbox{supp} \varrho_1 \subset B(0;1)$, define
$\varrho_{1,\varepsilon_1}(x)=\varepsilon_1^{-d}\varrho_1(x/
\varepsilon_1)$ and
$$
r_{\varepsilon_1}:= (b\cdot \nabla w) \ast\varrho_{1,\varepsilon_1} - b\cdot \nabla( w
\ast\varrho_{1,\varepsilon_1}).
$$
Let $Q \subset O \subset \subset \mathbb{R}^d$, with $Q$
compact and $O$ open, and let $L:= \|w\|_{L^\infty(O)}$. Then
\begin{eqnarray}\label{2.20}
\limsup_{\varepsilon_1\downarrow0}\int\limits_Q |r_{\varepsilon_1}|dx \leq L \int\limits_Q
\Lambda(M(x),\varrho_1(x))d|D^sb|(x) + L(d + I(\varrho_1)) |D^ab|(Q)
\end{eqnarray}
and
\begin{eqnarray}\label{2.21}
\limsup_{\varepsilon_1\downarrow0}\int\limits_Q |r_{\varepsilon_1}|dx \leq LI(\varrho_1)|D^s
b|(Q),
\end{eqnarray}
where $D^ab$ is the absolutely continuous part of the measure $Db$, and $D^sb$ is its singular part,
$$
\Lambda(M,\varrho_1):= \int\limits_{\mathbb{R}^d} | \langle Mz,\nabla\varrho_1(z)\rangle|
dz, \quad I(\varrho_1):= \int\limits_{\mathbb{R}^d}|z||\nabla\varrho_1(z)| dz,
$$
and $D^sb=M|D^sb|=\varsigma\otimes\vartheta|D^sb|$ with $\varsigma\bot\vartheta$ and $\varsigma, \vartheta \in S^{d-1}$.
 \vskip1mm\par
\textbf{(ii) (Anisotropic convolution kernels)} Given $\varsigma, \vartheta \in S^{d-1}$ with $\varsigma\bot\vartheta$  and given any
$\varepsilon_1> 0$ we can find a smooth, even convolution kernel
$\varrho_1$ with compact support such that
\begin{eqnarray}\label{2.22}
\int\limits_{\mathbb{R}^d}|\langle z,\varsigma\rangle||\langle\nabla \varrho_1(z),
\vartheta\rangle| dz <\varepsilon_1.
\end{eqnarray}
In particular $\Lambda(\varsigma\bot\vartheta, \varrho_1) < \varepsilon_1$.
\end{lemma}

\begin{lemma} [\cite{RY} Exercise 2.10, P31] \label{lem2.4} Let $\{X_t(x), x\in [0,1]^d, t\in [0,1]\}$ be a $\mathbb{R}^k$-valued continuous stochastic processes. Prove that if there exist three strictly positive constants $s,c, \varepsilon$ such that
\begin{eqnarray}\label{2.23}
\mathbb{E}[\sup_{0\leq t\leq 1}|X_t(x)-X_t(y)|^s]\leq c|x-y|^{d+\varepsilon},
\end{eqnarray}
then there is a modification $\tilde{X}$ of $X$ which is jointly continuous in $x$ and $t$ and is moreover H\"{o}lder continuous in $x$ of order $\beta$ for $\beta\in [0,\varepsilon/s)$ uniformly in $t$.

\vskip2mm\noindent
\textbf{Remark 2.2.} By following the proof of \cite[Theorem 2.1, page 25]{RY}, or the proof of \cite[Lemma 4.1]{TDW}, it is easy to prove above lemma, we omit the proof details. Moreover, from \cite[Theorem 4.1 (ii)]{TDW}, one has the following estimate
\begin{eqnarray}\label{2.24}
\mathbb{E}\sup_{0\leq t\leq 1}\|X_t(\cdot)\|^s_{\mathcal{C}_b^\beta(B_R)}<\infty.
\end{eqnarray}
\end{lemma}

\subsection{Main results}\label{sec2.3}

It is time for us to give a position to state our main results. The first result is concentrated on the existence and uniqueness of stochastic entropy solution for (\ref{1.6}).
\begin{theorem} \label{the2.1} \textbf{(Existence and uniqueness)} Let $\rho_0\in L^\infty\cap L^1(\mathbb{R}^d)$, $b\in {\mathcal C}_b^\alpha(\mathbb{R}^d;\mathbb{R}^d)$ with $\alpha\in (0,1)$, $f\in {\mathcal C}^1(\mathbb{R})$. We suppose in addition that
\begin{eqnarray}\label{2.25}
\mbox{div} b\in L^\infty(\mathbb{R}^d), \  f^\prime\in L^\infty(\mathbb{R}) \ \mbox{or} \ \mbox{div} b=0.
\end{eqnarray}
Then there exists a stochastic entropy solution to the Cauchy problem (\ref{1.6}). Moreover, if $b\in BV_{loc}(\mathbb{R}^d;\mathbb{R}^d)$, then the stochastic entropy solution is unique.
 \end{theorem}

\vskip2mm\noindent
\textbf{Remark 2.3.}  (i) When the noise vanishes, the general theory for entropy solutions was founded by Kru\u{z}kov \cite{Kru}. This theory is adequate for (\ref{1.2}), if we assume $b\in {\mathcal C}^1_b, f\in {\mathcal C}^1$ for uniqueness, and if we assume in addition that $b\in {\mathcal C}^3$ for existence. About twenty eighteen years later, using renormalized technique, DiPerna and Lions \cite{DL} also gained the existence and uniqueness for entropy solution but only assuming $b\in W^{1,1}_{loc}, \mbox{div} b\in L^\infty$, $|b(x)|/(1+|x|)\in L^1+L^\infty$ (the $W^{1,1}_{loc}$ regularity is sharp in the sense that divergence free vector-fields without integrable first derivatives exist counterexamples). Now, under the stochastic perturbation of a Brownian noise, when $\mbox{div} b\in L^\infty$, the boundedness and H\"{o}lder continuity of $b$ is sufficient to insure the existence of stochastic entropy solutions. Moreover, if $b\in BV_{loc}$ in addition, the uniqueness is also fulfilled. In this sense, we assert that the noise does have some regularizing effect.

(ii) When $f$ is linear, the Cauchy problem (\ref{1.6}) degenerate into the Cauchy problem (\ref{1.3}). For bounded vector field $b$ satisfying a globally H\"{o}lder continuous and an integrability condition on the divergence, in \cite{FGP}, Flandoli, Gubinelli and Priola showed the existence and uniqueness of $L^\infty$ solutions. As stated in \cite{FGP}, how to generalize the result to the nonlinear transport equation is still an open problem. The present result is an attempt to solve the problem and we extend the results given in \cite{AF,FF,FGP,LF} for the linear transport equations to the nonlinear transport equations.

Besides the existence and uniqueness of stochastic entropy solution, we are in a position to give a regularity result.

\begin{theorem}\label{the2.2} \textbf{(Regularity)} \ Let $b,f,\rho_0$ be described in Theorem \ref{the2.1}, which insure the existence and uniqueness for stochastic entropy solutions. If $\rho_0\in BV_x$ in addition, then the unique stochastic
entropy solution $\rho$ of (\ref{1.6}) is of class $L^1_{\omega}(L^\infty_{t}(BV_{loc}))$. Moreover, for almost all $\omega\in\Omega$, $\rho\in {\mathcal C}^\alpha([0,T];L^1(\mathbb{R}^d))$ for every $\alpha\in (0,1/2)$.
\end{theorem}

\vskip2mm\noindent
\textbf{Remark 2.4.}  For the deterministic equation, the existence of $BV$ solutions for (\ref{1.2}) in two space dimensions can be seen in \cite{Cag} if one supposes $b\in BV\cap L^\infty, \mbox{div} b=0$, $f\in{\mathcal C}^1$ and $f^\prime>0$. Now, with the noise, we prove that if $\mbox{div} b=0$, boundedness and H\"{o}lder continuity of $b$ and continuous differentiability of $f$ is enough to insure the existences of the $BV$ solution. And if one supposes $b\in BV_{loc}$ in addition, the uniqueness is conserved. Therefore, we extend the existing result but weaken the assumptions on $b$ and $f$.

To illustrate the regularizing effect, we present a counterexample for non-existence of such solution in the deterministic case. For simplicity, we suppose $f(\rho)=\rho$.

\begin{theorem} \label{the2.3} \textbf{(Non-existence)} Let $b$ and $\rho_0$ be stated in Theorem \ref{the2.2}. Consider the Cauchy problem
\begin{eqnarray}\label{2.26}
\left\{
  \begin{array}{ll}
\partial_t\rho(t,x)+b(x)\cdot\nabla \rho(t,x)=0, \ (t,x)\in (0,T)\times {\mathbb R}^d, \\
\rho(t,x)|_{t=0}=\rho_0(x), \  x\in{\mathbb R}^d.
  \end{array}
\right.
\end{eqnarray}
Then there exists a unique weak solution. However, if $d\geq 2$, one can choose proper functions $\rho_0$ and $b$ such
that the unique weak solution $\rho(t,x)$ does not lie in
$L^\infty_t(BV_{loc})$. Here $\rho$ is said to be a weak
solution of (\ref{2.26}), if it lies in $L^\infty_{t,x}\cap {\mathcal C}_t(L^1_x)$ and meets (\ref{2.26}) in the sense of distributions.
 \end{theorem}

\vskip2mm\noindent
\textbf{Remark 2.5.} The phenomenon that a PDE of fluid dynamics that while existence fails, then a random perturbation rends the equation well-posed,  occurs not only in a transport equations (regarded as a first order hyperbolic equation) but it also happens in a wave equation (regarded as a second order hyperbolic equation).  Indeed, when the transport equation in (\ref{2.26}) is replaced by a wave equation and $\mathbb{R}^d$ is replaced by a three dimensional compact Riemann manifold $M$, with an additional nonhomogeneous term, it becomes
\begin{eqnarray}\label{2.27}
\left\{
  \begin{array}{ll}
\partial^2_t\rho(t,x)-\Delta\rho(t,x)+\rho^3=0, \ \ t>0, \ x\in M,  \\
\rho(t,x)|_{t=0}=\rho_0(x), \ x\in M.
  \end{array}
\right.
\end{eqnarray}
When the initial data in $H^s(M), s < 1/2$, this problem is supercritical and can be shown to be strongly ill-posed (in the Hadamard sense).  However, after a suitable randomization, Bouq and Tzvetkov \cite{BT1} showed the local existence of strong solutions for (\ref{2.27}) for a large set of initial data in $H^s(M), s \geq1/4$. For more details in this direction, one also consults  to \cite{BT2}.

\section{Existence}\label{sec3}
\setcounter{equation}{0}
In this section, we prove the first part of Theorem \ref{the2.1}. That is the following
\begin{theorem} \label{the3.1} Let $\rho_0\in L^\infty\cap L^1(\mathbb{R}^d)$, $b\in {\mathcal C}_b^\alpha(\mathbb{R}^d;\mathbb{R}^d)$ with $\alpha\in (0,1)$, $f\in {\mathcal C}^1(\mathbb{R})$, that (\ref{2.25}) holds. Then there exists a stochastic entropy solution to the Cauchy problem (\ref{1.6}).
\end{theorem}
\vskip2mm\noindent
\textbf{Proof.}  Since the proof for $\mbox{div} b=0$ is similar to the case of $\mbox{div} b\in L^\infty(\mathbb{R}^d),  f^\prime\in L^\infty(\mathbb{R})$ (in fact, it is easier), we only give the details for $\mbox{div} b\in L^\infty$.

Inspired by Lemma \ref{lem2.1}, it suffices to establish the existence of stochastic weak solutions for (\ref{2.2}). We found the existence of stochastic weak solutions for (\ref{2.2}) by using the stochastic BGK approximation, i.e. we regard (\ref{2.2}) as the $\varepsilon_1\downarrow 0$ limit of the
integro-differential equation
\begin{eqnarray}\label{3.1}
\left\{
  \begin{array}{ll}
\partial_tu_{\varepsilon_1}(t,x,v)+f^\prime(v)b(x)\cdot\nabla_x
u_{\varepsilon_1}(t,x,v)
+\partial_{x_i}u_{\varepsilon_1}\circ \dot{B}_i(t)= \frac{1}{\varepsilon_1}
{\big[} \chi_{\rho_{\varepsilon_1}} -u_{\varepsilon_1}{\big]}, \\
u_{\varepsilon_1}(t,x,v)|_{t=0}=\chi_{\rho_0(x)}(v),  \\ \rho_\varepsilon=\rho_\varepsilon(t,x)=\int_\mathbb{R} u_{\varepsilon_1}(t,x,v)dv.
  \end{array}
\right.
\end{eqnarray}
Here $u_{\varepsilon_1}$ is said to be a stochastic weak solution of (\ref{3.1}) if it is a predictable random field and lies in $L^\infty_{\omega,t,x}\cap {\mathcal C}_t(L^1_{\omega,x,v})$. For every $\varphi\in {\mathcal D}(\mathbb{R}^d)$, every $\psi\in {\mathcal D}(\mathbb{R})$, every $t\in [0,T)$ and for almost all $\omega\in\Omega$,
\begin{eqnarray}\label{3.2}
&&\int\limits_{\mathbb{R}^{d+1}_{x,v}}\varphi(x)\psi(v)u_{\varepsilon_1}(t,x,v)dxdv\cr\cr&=& \int\limits_{\mathbb{R}^{d+1}_{x,v}}\varphi(x)\psi(v)\chi_{\rho_0(x)}(v)dxdv+
\int\limits^t_0\int\limits_{\mathbb{R}^{d+1}_{x,v}}f^\prime(v)\mbox{div}(b(x)\varphi(x))\psi(v)
u_{\varepsilon_1}(s,x,v)dxdvds\cr\cr&&+\int\limits^t_0\circ
dB_i(s)\int\limits_{\mathbb{R}^{d+1}_{x,v}}\psi(v)\partial_{x_i}
\varphi(x)u_{\varepsilon_1}dxdv+\frac{1}{\varepsilon_1}\int\limits^t_0\int\limits_{\mathbb{R}^{d+1}_{x,v}}
\psi(v)\varphi(x)[\chi_{\rho_{\varepsilon_1}}-u_{\varepsilon_1}]dxdvds.
\end{eqnarray}

The proof details can be given into five steps.
 \vskip2mm\noindent
\textbf{Step 1.} The existence and uniqueness of the strong solution for SDE:
\begin{eqnarray}\label{3.3}
dX(s,t)=f^\prime(v)b(X(s,t))dt+dB(t),  \ t\in(s,T], \ X(s,t)|_{t=s}=x\in \mathbb{R}^d,
\end{eqnarray}
where $v\in\mathbb{R}$ is a parameter.

To prove the result, we use It\^{o}-Tanack's trick (see see
\cite{FGP} for more details). Initially, we argue the following backward second order parabolic equation:
\begin{eqnarray}\label{3.4}
\left\{
\begin{array}{ll}
\partial_{t}U(t,x) +\frac{1}{2}\Delta U(t,x)+f^\prime(v)b(x)\cdot \nabla U(t,x)=
\lambda U(t,x)-f^\prime(v)b(x), \ (t,x)\in (s,T)\times {\mathbb R}^d, \\
U(T,x)=0, \  x\in{\mathbb R}^d.
  \end{array}
\right.
\end{eqnarray}
Due to the hypotheses in Theorem \ref{the3.1} on $f$ and $b$, with the help of classical parabolic theory (see \cite{Kry}), for every $v\in \mathbb R$, there is a unique classical solution $U$ to (\ref{3.4}). Moreover, by a direct computation, we obtain that $U\in \mathcal{C}_b([s,T]\times{\mathbb R}_v;\mathcal{C}_b^{2+\alpha}({\mathbb R}^d;{\mathbb R}^d))\cap \mathcal{C}_b({\mathbb R}_v;\mathcal{C}^1([s,T];\mathcal{C}_b^\alpha({\mathbb R}^d;{\mathbb R}^d)))$, and
\begin{eqnarray}\label{3.5}
\|\nabla_xU\|_{L^\infty((s,T)\times{\mathbb R}_v\times {\mathbb R}^d)} \rightarrow 0, \ as \ \lambda\rightarrow \infty.
\end{eqnarray}

For every $v\in{\mathbb R}$, we set $U_v(t,x)=U(t,x,v)$ and $\gamma_v(t,x)=x+U_v(t,x)$, then in view of (\ref{3.5}) and recall the classical Hadamard theorem (see \cite[page 32]{Pro}), for a given and big enough real number $\lambda$, for every $t\in[s,T]$, and every $v\in {\mathbb R}$, $\gamma_v(t,\cdot)$ forms a non-singular diffeomorphism of $\mathcal{C}^2$. Consider the following SDE:
\begin{eqnarray}\label{3.6}
\left\{
\begin{array}{ll}
d Y^v(s,t)=\lambda U_v(t,\gamma^{-1}_v(t,Y^v(s,t)))dt+ [I+\nabla_y
U_v(t,\gamma^{-1}_v(t,Y^v(s,t)))] dB(t), \ t\in(s,T], \\ Y^v(s,t)|_{t=s}=y,
 \end{array}
\right.
\end{eqnarray}
where $\gamma^{-1}_v(t,x)$ is its inverse of
the mapping $x\mapsto\gamma_v(t,x)$. Since $\gamma_v$ and $\gamma^{-1}_v$ have a bounded first and second spatial derives, with the help of classical result (see \cite{Kun}), for every $v\in {\mathbb R}$, there is a unique strong solution $Y^v(s,t,y)$ of (\ref{3.6}), which forms stochastic homeomorphism flow of class ${\mathcal C}^{1,\alpha^\prime}$ ($\alpha^\prime<\alpha$). Moreover, for every $x,y\in\mathbb{R}^d$, by a direct computations (see \cite{WDGL} for example), for every $p\geq1$, we have
\begin{eqnarray}\label{3.7}
\sup_{v\in\mathbb{R}}{\mathbb E}|Y^v(s,t,x)-Y^v(s,t,y)|^p\leq C|x-y|^p,
\end{eqnarray}
and
\begin{eqnarray}\label{3.8}
\sup_{v\in\mathbb{R}}{\mathbb E}|\nabla Y^v(s,t,x)-\nabla Y^v(s,t,y)|^p\leq C|x-y|^{\alpha p}.
\end{eqnarray}
From (\ref{3.7}) and (\ref{3.8}), thank to Lemma \ref{lem2.4} and Remark 2.2, so (\ref{2.24}) holds true for the present stochastic process, i.e.
$Y^\cdot(s,\cdot,\cdot)\in {\mathcal C}_b(\mathbb{R}_v;L^p(\Omega;\mathcal{C}_b([s,T];
\mathcal{C}_b^{1+\alpha^\prime}(B_R))))$ for every $p,R\geq 1$.

On the other hand, by the relationship $X^v(s,t)=\gamma^{-1}_v\circ Y^v(s,t)$, for every $v\in \mathbb{R}$, there is a unique strong solution $X^v(s,t,x)$ of (\ref{3.3}), which forms stochastic homeomorphism flow of class  ${\mathcal C}^{1,\alpha^\prime}$ ($\alpha^\prime<\alpha$). Besides, $X^\cdot(s,\cdot,\cdot)\in {\mathcal C}_b(\mathbb{R}_v;L^p(\Omega;\mathcal{C}_b([s,T];
\mathcal{C}_b^{1+\alpha^\prime}(B_R))))$ for every $p,R\geq 1$, so we show that $X^v(s,t,x)$ is measurable in $v$. Further, by the Liouville theorem, we have the following Euler's identity:
\begin{eqnarray}\label{3.9}
|\nabla_xX^v(s,t,x)|=\exp(f^\prime(v)\int\limits^t_s\mbox{div} b(X^v(s,r,x))dr=:J(s,t,X),
\end{eqnarray}
where the $|\cdot|$ in the left hand of (\ref{3.3}) denotes the  determinant of the matrix. Moreover, by the representation of (\ref{3.3}), for almost all $\omega\in\Omega$, $X^\cdot(s,t,x,\omega)\in
\mathcal{C}_b({\mathbb R}_v)$.

\vskip2mm\noindent
\textbf{Step 2.} The existence of stochastic weak solution for (\ref{3.1}).
\vskip1.5mm\par
If $u_{\varepsilon_1}$ is a classical solution of (\ref{3.1}) (i.e. continuously differentiable in $t$ and $x$),  along the direction (\ref{3.3}), then
\begin{eqnarray}\label{3.10}
u_{\varepsilon_1}(t,X^v(0,t,x),v)=\frac{1}{\varepsilon_1}\int\limits^t_0e^{\frac{s-t}{\varepsilon_1}}
\chi_{\rho_{\varepsilon_1}(s,X^v(0,s,x))}(v)ds+
e^{-\frac{t}{\varepsilon_1}}\chi_{\rho_0(x)}(v).
\end{eqnarray}
Thus
\begin{eqnarray}\label{3.11}
u_{\varepsilon_1}(t,x,v)=\frac{1}{\varepsilon_1}\int\limits^t_0e^{\frac{s-t}
{\varepsilon_1}}
\chi_{\rho_{\varepsilon_1}(s,X_{t,s}^v(x))}(v)ds+
e^{-\frac{t}{\varepsilon_1}}\chi_{\rho_0(X^v_{t,0}(x))}(v),
\end{eqnarray}
where
$X_{t,s}^v(x)=(X_{s,t}^v(x))^{-1}=(X^v(s,t,x))^{-1}$.

Now, we will show that if $u_{\varepsilon_1}\in L^\infty_{\omega,t,x}\cap {\mathcal C}_t(L^1_{\omega,x,v})$ is a predictable random field, which satisfies (\ref{3.11}), then $u_{\varepsilon_1}$ is a stochastic weak solution of (\ref{3.1}).

Let $u_{\varepsilon_1}$ given by (\ref{3.11}) and $X^v(s,t,x)$ meets (\ref{3.3}), then for every $\varphi\in {\mathcal D}(\mathbb{R}^d)$ and for every $\psi\in {\mathcal D}(\mathbb{R})$,
\begin{eqnarray}\label{3.12}
&&\int\limits_{\mathbb{R}^{d+1}_{x,v}}\varphi(x)\psi(v)u_{\varepsilon_1}(t,x,v)dxdv
\cr\cr&=& \frac{1}{\varepsilon_1}\int\limits^t_0e^{\frac{s-t}
{\varepsilon_1}}\int\limits_{\mathbb{R}^{d+1}_{x,v}}\varphi(x)\psi(v)
\chi_{\rho_{\varepsilon_1}(s,X_{t,s}^v)}(v)dxdvds
+
e^{-\frac{t}{\varepsilon_1}} \int\limits_{\mathbb{R}^{d+1}_{x,v}}\varphi(x)\psi(v)\chi_{\rho_0(X^v_{t,0})}(v)dxdv
\cr\cr&=&\frac{1}{\varepsilon_1}\int\limits^t_0
\int\limits_{\mathbb{R}^{d+1}_{x,v}}I(\varphi,J,X^v_{s,t})\psi(v)
\chi_{\rho_{\varepsilon_1}(s,x)}(v)dxdvds
\cr\cr&&+
\int\limits_{\mathbb{R}^{d+1}_{x,v}}I(\varphi,J,X^v_{0,t})\psi(v)\chi_{\rho_0(x)}(v)dxdv, \quad {\mathbb P}-a.s.,
\end{eqnarray}
where for every $0\leq r\leq t\leq T$,
$$
I(\varphi,J,X^v_{r,t})=e^{\frac{r-t}{\varepsilon_1}}\varphi(X^v_{r,t}(x))J(r,t,X).
$$
Using It\^{o}'s formula, then we have
\begin{eqnarray}\label{3.13}
I(\varphi,J,X^v_{r,t})&=&\varphi(x)+
\int\limits_r^tI(\varphi,J,X^v_{r,\tau})\Big[f^\prime(v)\mbox{div} b(X^v_{r,\tau}(x))-\frac{1}{\varepsilon_1}\Big]d\tau
\cr\cr&&+\int\limits_r^te^{\frac{r-\tau}{\varepsilon_1}}J(r,\tau,X)
f^\prime(v)b(X^v_{r,\tau}(x))\cdot\nabla_x\varphi(X^v_{r,\tau}(x))d\tau
\cr\cr&&+\int\limits_r^te^{\frac{r-\tau}{\varepsilon_1}}J(r,\tau,X)
\partial_{x_i}\varphi(X^v_{r,\tau}(x))\circ dB_i(\tau).
\end{eqnarray}
Combining (\ref{3.12}) and (\ref{3.13}), one ends up with
\begin{eqnarray}\label{3.14}
&&\int\limits_{\mathbb{R}^{d+1}_{x,v}}\varphi(x)\psi(v)u_{\varepsilon_1}(t,x,v)dxdv
\cr\cr&=&\frac{1}{\varepsilon_1}\int\limits^t_0
\int\limits_{\mathbb{R}^{d+1}_{x,v}}\varphi(x)\psi(v)
\chi_{\rho_{\varepsilon_1}(\tau,x)}(v)dxdvd\tau
\cr\cr&&+
\int\limits^t_0
\int\limits_{\mathbb{R}^{d+1}_{x,v}}\psi(v)\Big[ f^\prime(v)\mbox{div} (b(x)\varphi(x))-\frac{1}{\varepsilon_1}\varphi(x)\Big]\Big[
\frac{1}{\varepsilon_1}\int\limits^\tau_0 e^{\frac{s-\tau}{\varepsilon_1}}\chi_{\rho_{\varepsilon_1}(s,X^v_{\tau,s})}(v)
ds\Big]dxdvd\tau
\cr\cr&&
+\int\limits^t_0\circ dB_i(\tau)
\int\limits_{\mathbb{R}^{d+1}_{x,v}}
\Big[\frac{1}{\varepsilon_1}\int\limits_0^\tau e^{\frac{s-\tau}{\varepsilon_1}}
\chi_{\rho_{\varepsilon_1}(s,X^v_{\tau,s})}(v)
ds\Big]
\partial_{x_i}\varphi(x)
\psi(v)dxdv
\cr\cr&&+
\int\limits_{\mathbb{R}^{d+1}_{x,v}}\varphi(x)\psi(v)
\chi_{\rho_{0}(x)}(v)dxdv
\cr\cr&&+\int\limits_0^t
\int\limits_{\mathbb{R}^{d+1}_{x,v}}
 e^{-\frac{\tau}{\varepsilon_1}}\psi(v)\Big[f^\prime(v)\mbox{div} (b(x)\varphi(x))-\frac{1}{\varepsilon_1}\varphi(x)\Big]
\chi_{\rho_0(X^v_{\tau,0})}(v)dxdvd\tau
\cr\cr&&+\int\limits_0^t\circ dB_i(\tau)
\int\limits_{\mathbb{R}^{d+1}_{x,v}}e^{-\frac{\tau}{\varepsilon_1}}
\partial_{x_i}\varphi(x)
\psi(v)\chi_{\rho_0(X^v_{\tau,0})}(v)dxdv
\cr\cr&=&
\int\limits_{\mathbb{R}^{d+1}_{x,v}}\varphi(x)\psi(v)
\chi_{\rho_{0}(x)}(v)dxdv\cr\cr&&+
\int\limits^t_0
\int\limits_{\mathbb{R}^{d+1}_{x,v}}\psi(v)f^\prime(v)\mbox{div} (b(x)\varphi(x))\Big[
\frac{1}{\varepsilon_1}\int\limits^\tau_0 e^{\frac{s-\tau}{\varepsilon_1}}\chi_{\rho_{\varepsilon_1}(s,X^v_{\tau,s})}(v)
ds+e^{-\frac{\tau}{\varepsilon_1}}
\chi_{\rho_0(X^v_{\tau,0})}(v)\Big]dxdvd\tau
\cr\cr&&
+\int\limits^t_0\circ dB_i(\tau)
\int\limits_{\mathbb{R}^{d+1}_{x,v}}
\Big[\frac{1}{\varepsilon_1}\int\limits_0^\tau e^{\frac{s-\tau}{\varepsilon_1}}
\chi_{\rho_{\varepsilon_1}(s,X^v_{\tau,s})}(v)
ds+e^{-\frac{\tau}{\varepsilon_1}}\chi_{\rho_0(X^v_{\tau,0})}(v)\Big]
\partial_{x_i}\varphi(x)
\psi(v) dxdv
\cr\cr&&+
\frac{1}{\varepsilon_1}\int\limits^t_0
\int\limits_{\mathbb{R}^{d+1}_{x,v}}\varphi(x)\psi(v)
\chi_{\rho_{\varepsilon_1}(\tau,x)}(v)dxdvd\tau
\cr\cr&&
-\frac{1}{\varepsilon_1}\int\limits^t_0
\int\limits_{\mathbb{R}^{d+1}_{x,v}}\varphi(x)\psi(v)\Big[
\frac{1}{\varepsilon_1}\int\limits^\tau_0 e^{\frac{s-\tau}{\varepsilon_1}}\chi_{\rho_{\varepsilon_1}(s,X^v_{\tau,s})}(v)
ds+e^{-\frac{\tau}{\varepsilon_1}}
\chi_{\rho_0(X^v_{\tau,0})}(v)\Big]dxdvd\tau.
\end{eqnarray}
By virtue of (\ref{3.11}), from (\ref{3.14}), one achieves identity (\ref{3.2}), so $u_{\varepsilon_1}$ is a stochastic weak solution of (\ref{3.1}). It remains to verify the existence of stochastic weak solution for (\ref{3.11}).

For every $u\in L^\infty_{\omega,t,x}\cap {\mathcal C}_t(L^1_{\omega,x,v})$, which is predictable, defining
$S_{\varepsilon_1}$ by:
\begin{eqnarray}\label{3.15}
(S_{\varepsilon_1} u)(t,x,v)=\frac{1}{\varepsilon_1}\int\limits^t_0e^{\frac{s-t}
{\varepsilon_1}} \chi_{\rho^u(s,X^v_{t,s}(x))}(v)ds+
e^{-\frac{t}{\varepsilon_1}}\chi_{\rho_0^u(X^v_{t,0}(x))}(v),
\end{eqnarray}
with
$$
\rho^u(t,x)=\int\limits_{\mathbb{R}}u(t,x,v)dv, \ \ \ \rho_0^u(x)=\int\limits_{{\mathbb
R}}u_0(x,v)dv=\rho_0(x).
$$
Then the random field $S_{\varepsilon_1}u$ is predictable. Now, we collate that (\ref{3.15}) is well-defined. Indeed,
\begin{eqnarray}\label{3.16}
\|S_{\varepsilon_1} u\|_{L^\infty_{\omega,t,x,v}}\leq
1.
\end{eqnarray}
With the help of $\mbox{div} b\in L^\infty$ and $f^\prime\in L^\infty$,
\begin{eqnarray}\label{3.17}
&&\|S_{\varepsilon_1} u\|_{{\mathcal C}_t(L^1_{\omega,x,v})}
\cr\cr&\leq&\sup_{0\leq t\leq T}{\mathbb E}\Big[
\frac{1}{\varepsilon_1}\int\limits^t_0e^{\frac{s-t}
{\varepsilon_1}}ds\int\limits_{\mathbb{R}^{d+1}_{x,v}}|\chi_{\rho^u(s,X^v_{t,s}(x))}(v)|dxdv +
e^{-\frac{t}{\varepsilon_1}}\int\limits_{\mathbb{R}^{d+1}_{x,v}}
|\chi_{\rho_0^u(X_{t,0}^v(x))}(v)|dxdv\Big] \cr \cr &\leq&
\sup_{0\leq t\leq T}{\mathbb E}\Big[
\frac{1}{\varepsilon_1}\int\limits^t_0e^{\frac{s-t}
{\varepsilon_1}}ds\int\limits_{\mathbb{R}^{d+1}_{x,v}}|\chi_{\rho^u(s,x)}(v)|
\exp(f^\prime(v)\int\limits^t_s\mbox{div} b(X^v_{s,r}(x))dr) dxdv \cr \cr && +
e^{-\frac{t}{\varepsilon}}\int\limits_{\mathbb{R}^{d+1}_{x,v}}
|\chi_{\rho_0^u(x)}(v)|\exp(f^\prime(v)\int\limits^t_0
\mbox{div} b(X^v_{0,r}(x))dr)dxdv\Big] \cr \cr &\leq&
\exp(T\|f^\prime\|_{L^\infty_v}\|\mbox{div} b\|_{L^\infty_x})\Big[(1-e^{-\frac{T}{\varepsilon}})\|u\|_{{\mathcal C}_t(L^1_{\omega,x,v})}+
\|u_0\|_{L^1_{x,v}}\Big].
\end{eqnarray}

By (\ref{3.16}) and (\ref{3.17}), thus (\ref{3.15}) is well-defined. Moreover, for every $f,g\in L^\infty_{\omega,t,x}\cap {\mathcal C}_t(L^1_{\omega,x,v})$, and every $0<T_1\leq T$, we have the following estimate
\begin{eqnarray}\label{3.18}
 && \|S_{\varepsilon_1} f-S_{\varepsilon_1} g\|_{{\mathcal C}([0,T_1];
L^1(\Omega\times\mathbb{R}^{d+1}_{x,v}))}
 \cr \cr
&\leq& \sup_{0\leq t\leq T_1}{\mathbb E}\Big[
\frac{1}{\varepsilon_1}\int\limits^t_0e^{\frac{s-t}
{\varepsilon_1}}ds\int\limits_{\mathbb{R}^{d+1}_{x,v}}
|\chi_{\rho^f(s,X^v_{t,s}(x))}(v)-
\chi_{\rho^g(s,X^v_{t,s}(x))}(v)| dxdv
 \cr \cr &&+
e^{-\frac{t}{\varepsilon_1}}\int\limits_{\mathbb{R}^{d+1}_{x,v}}
|\chi_{\rho_0^f(X_{t,0}^v(v))}(v)-
\chi_{\rho_0^g(X_{t,0}^v(v))}(v)|dxdv\Big] \cr \cr &=&
\sup_{0\leq t\leq T_1}{\mathbb E}\Big[
\frac{1}{\varepsilon_1}\int\limits^t_0e^{\frac{s-t}
{\varepsilon_1}}ds\int\limits_{\mathbb{R}^{d+1}_{x,v}}
|\chi_{\rho^f(s,x)}(v)-\chi_{\rho^g(s,x)}(v)|
\exp(f^\prime(v)\int\limits^t_s\mbox{div} b(X^v_{s,r})dr) dxdv \cr \cr && + \
e^{-\frac{t}{\varepsilon_1}}\int\limits_{\mathbb{R}^{d+1}_{x,v}}
|\chi_{\rho_0^f(x)}(v)-\chi_{\rho_0^g(x)}(v) |
\exp(f^\prime(v)\int\limits^t_0\mbox{div} b(X_{0,r}^v(x))dr)dxdv\Big] \cr \cr
&\leq& \exp(T_1\|f^\prime\|_{L^\infty_v}\|\mbox{div} b\|_{L^\infty_x})\Big[
(1-e^{-\frac{T_1}{\varepsilon_1}})\|f-g\|_{{\mathcal C}([0,T_1];
L^1(\Omega\times\mathbb{R}^{d+1}_{x,v}))}+\|f_0-g_0\|_{L^1_{x,v}}\Big].
\end{eqnarray}
In particular, if $f_0=g_0=\chi_{\rho_0}$, from (\ref{3.12}), we get
$$
\|S_{\varepsilon_1} f-S_{\varepsilon_1} g\|_{{\mathcal C}([0,T_1];
L^1(\Omega\times\mathbb{R}^{d+1}_{x,v}))}\leq \exp(T_1\|f^\prime\|_{L^\infty_v}\|\mbox{div} b\|_{L^\infty_x})
(1-e^{-\frac{T_1}{\varepsilon_1}})\|f-g\|_{{\mathcal C}([0,T_1];
L^1(\Omega\times\mathbb{R}^{d+1}_{x,v}))}.
$$

Since $0<T_1\leq T$ is arbitrary, we can select $T_1$ so small that
$\exp(T_1\|f^\prime\|_{L^\infty_v}\|\mbox{div} b\|_{L^\infty_x})
(1-e^{-\frac{T_1}{\varepsilon_1}})<1$. We apply the Banach fixed
point theorem to find a unique $u_{\varepsilon_1} \in
{\mathcal C}([0,T_1];L^1(\Omega\times\mathbb{R}^{d+1}_{x,v}))$ solving the Cauchy problem (\ref{3.11}). By (\ref{3.16}), $u_{\varepsilon_1} \in
L^\infty([0,T_1];L^\infty(\Omega\times\mathbb{R}^{d+1}_{x,v}))$, so
$u_{\varepsilon_1}(T_1)\in L^1(\Omega\times\mathbb{R}^{d+1}_{x,v})\cap
L^\infty(\Omega\times\mathbb{R}^{d+1}_{x,v})$. We then repeat the argument above to extend our solution to the time interval $[T_1,2T_1]$. Continuing, after finitely many steps we construct a solution existing on the interval $[0,T]$. From this, we demonstrate that there
exists a unique $u_{\varepsilon_1} \in
{\mathcal C}([0,T];L^1(\Omega\times\mathbb{R}^{d+1}_{x,v}))\cap
L^\infty([0,T]\times\Omega\times\mathbb{R}^{d+1}_{x,v})$ solving the
Cauchy problem (\ref{3.11}).

Thanks to (\ref{3.11}) and (\ref{3.17}), one gains the following estimate
\begin{eqnarray*}
\|u_{\varepsilon_1}(t)\|_{L^1_{x,v}}
&\leq&
\int\limits_{\mathbb{R}^{d+1}_{x,v}}  \frac{1}{\varepsilon_1}\int\limits^t_0e^{\frac{s-t}
{\varepsilon_1}}e^{(t-s)C_0}|
\chi_{\rho_{\varepsilon_1}(s,x)}(v)|dsdxdv+
\int\limits_{\mathbb{R}^{d+1}_{x,v}}e^{-\frac{t}{\varepsilon_1}}e^{C_0t}
|\chi_{\rho_0(x)}(v)|dxdv\cr\cr &\leq&
\frac{1}{\varepsilon_1}\int\limits^t_0e^{\frac{s-t}
{\varepsilon_1}}e^{(t-s)C_0}\|u_{\varepsilon_1}(s)\|_{L^1_{x,v}}ds
+e^{-\frac{t}{\varepsilon_1}}e^{C_0t}
\|\rho_0\|_{L^1_x}, \quad {\mathbb P}-a.s.,
\end{eqnarray*}
which yields that
\begin{eqnarray*}
U(t)\leq (1-e^{-\frac{t}{\varepsilon_1}})\max_{0\leq s\leq t}U(s)+e^{-\frac{t}{\varepsilon_1}}U(0), \quad {\mathbb P}-a.s.,
\end{eqnarray*}
where $U(t)=e^{-C_0t}\|u_{\varepsilon_1}(t)\|_{L^1_{x,v}}$, $C_0=\|f^\prime\|_{L^\infty_v}\|\mbox{div} b\|_{L^\infty_x}$. It follows that
\begin{eqnarray}\label{3.19}
\|u_{\varepsilon_1}(t)\|_{L^1_{x,v}}
\leq e^{C_0T}
\|\rho_0\|_{L^1_x}, \quad {\mathbb P}-a.s..
\end{eqnarray}

\vskip2mm\noindent
\textbf{Step 3.} $u_{\varepsilon_1} \in L^\infty_{\omega,t,x}(L^1_v)$ and
\begin{eqnarray}\label{3.20}
\|\rho_{\varepsilon_1}(t)\|_{L^\infty_{\omega,t,x}}\leq \|u_{\varepsilon_1}(t)\|_{L^\infty_{\omega,t,x}(L^1_v)}\leq\|\rho_0\|_{L^\infty_x}.
\end{eqnarray}

Obviously, the first inequality is natural. It is sufficient to show the second inequality. In fact, from (\ref{3.11}), for almost all $\omega\in\Omega$,
\begin{eqnarray}\label{3.21}
\|u_{\varepsilon_1}(t,\cdot,v)\|_{L^\infty_x}
\leq
\frac{1}{\varepsilon_1}\int\limits^t_0e^{\frac{s-t}
{\varepsilon_1}}\|\chi_{\rho_{\varepsilon_1}(s,\cdot)}(v)\|_{L^\infty_x}ds
+
e^{-\frac{t}{\varepsilon_1}}\|\chi_{\rho_0}(v)\|_{L^\infty_x}.
\end{eqnarray}

If one integrates the (\ref{3.21}) with respect to the variable $v$,
\begin{eqnarray}\label{3.22}
&&\int\limits_{\mathbb{R}}\|u_{\varepsilon_1}(t,\cdot,v)\|_{L^\infty_x}dv
\cr\cr&\leq&\frac{1}{\varepsilon_1}\int\limits^t_0e^{\frac{s-t}
{\varepsilon_1}} \int\limits_{\mathbb{R}}\|\chi_{\rho_{\varepsilon_1}(s,\cdot)}(v)\|_{L^\infty_x}dvds+
e^{-\frac{t}{\varepsilon_1}}\int\limits_{\mathbb{R}}\|\chi_{\rho_0}(v)\|_{L^\infty_x}dv
\cr\cr&=&\frac{1}{\varepsilon_1}\int\limits^t_0e^{\frac{s-t}
{\varepsilon_1}} \int\limits_{\mathbb{R}}\max\{\|\textbf{1}_{(0,\rho_{\varepsilon_1}^+)}(v)\|_{L^\infty_x},  \|\textbf{1}_{(-\rho_{\varepsilon_1}^-,0)}(v)\|_{L^\infty_x}\}dvds
\cr\cr&&+e^{-\frac{t}{\varepsilon_1}}\int\limits_{\mathbb{R}}
\|\chi_{\rho_0}(v)\|_{L^\infty_x}dv, \quad {\mathbb P}-a.s.,
\end{eqnarray}
where in the last line in (\ref{3.22}), we have used the definition of the function $\chi_\cdot(\cdot)$ (see Section 2, after (\ref{2.2})).

Observing the fact that $\|\textbf{1}_{(0,\rho_{\varepsilon_1}^+)}(\cdot)\|_{L^\infty_x}  \|\textbf{1}_{(-\rho_{\varepsilon_1}^-,0)}(\cdot)\|_{L^\infty_x}=0$, from (\ref{3.22}), then
\begin{eqnarray}\label{3.23}
&&\int\limits_{\mathbb{R}}\|u_{\varepsilon_1}(t,\cdot,v)\|_{L^\infty_x}dv
\cr\cr&\leq&\frac{1}{\varepsilon_1}\int\limits^t_0e^{\frac{s-t}
{\varepsilon_1}}\max\{\|\rho_{\varepsilon_1}^+\|_{L^\infty_x}, \|\rho_{\varepsilon_1}^-\|_{L^\infty_x}\}ds
+e^{-\frac{t}{\varepsilon_1}}\int\limits_{\mathbb{R}}\|\chi_{\rho_0}(v)\|_{L^\infty_x}dv
\cr\cr&=&\frac{1}{\varepsilon_1}\int\limits^t_0e^{\frac{s-t}
{\varepsilon_1}}{\Big \|}\int\limits_{\mathbb{R}}u_{\varepsilon_1}(t,x,v)dv{\Big \|}_{L^\infty_x}ds
+e^{-\frac{t}{\varepsilon_1}}\int\limits_{\mathbb{R}}\|\chi_{\rho_0}(v)\|_{L^\infty_x}dv
\cr\cr&\leq&(1-e^{-\frac{t}{\varepsilon_1}})
\max_{0\leq s\leq t} \int\limits_{\mathbb{R}}\|u_{\varepsilon_1}(s,\cdot,v)\|_{L^\infty_x}dv+
e^{-\frac{t}{\varepsilon_1}}\int\limits_{\mathbb{R}}\|\chi_{\rho_0}(v)\|_{L^\infty_x}dv,
\quad {\mathbb P}-a.s.,
\end{eqnarray}
where in the last inequality, we used the Minkowski inequality
$$
\|u_{\varepsilon_1}(t,\cdot,\cdot)\|_{L^\infty_x(L^1_v)}\leq \|u_{\varepsilon_1}(t,\cdot,\cdot)\|_{L^1_v(L^\infty_x)}.
$$

So we find that the function $J(t)$,
$$
J(t)=\int\limits_{\mathbb{R}}\|u_{\varepsilon_1}(t,\cdot,v)\|_{L^\infty_x}dv,
$$
satisfies
$$
J(t)\leq (1-e^{-\frac{t}{\varepsilon_1}})\max_{0\leq s\leq t}J(s)+e^{-\frac{t}{\varepsilon_1}}J(0), \quad {\mathbb P}-a.s.,
$$
which implies
\begin{eqnarray*}
J(t)\leq J(0).
\end{eqnarray*}

With the help of the Minkowski inequality and the boundedness of $\rho_0$, from (\ref{3.23}),
\begin{eqnarray*}
\|u_{\varepsilon_1}(t)\|_{L^\infty_{\omega,x}(L^1_v)}\leq \|u_{\varepsilon_1}(t)\|_{L^\infty_\omega(L^1_v(L^\infty_x))}\leq \|\chi_{\rho_0}\|_{L^1_v(L^\infty_x)}\leq \|\rho_0\|_{L^\infty_x}.
\end{eqnarray*}

\vskip2mm\noindent
\textbf{Step 4.} The $\varepsilon_1\downarrow 0$ limits for approximate solutions.

In view of  (\ref{3.16}) and (\ref{3.20}), there exist two subsequences (denoted them by themselves for simplicity) $\{u_{\varepsilon_1}\}$ and $\{\rho_{\varepsilon_1}\}$,
such that
\begin{eqnarray}\label{3.24}
u_{\varepsilon_1}\stackrel{w*}{\longrightarrow} u \ \ \mbox{in} \ \ L^\infty_{\omega,t,x,v}, \ \
\rho_{\varepsilon_1}\stackrel{w*}{\longrightarrow} \rho \ \ \mbox{in} \ \  L^\infty_{\omega,t,x},  \ \  \mbox{as} \ \ \varepsilon_1\rightarrow 0.
\end{eqnarray}
Since the space of predictable process is weakly-closed, $u$ and $\rho$ are predictable. On the other hand, by (\ref{3.19}) and the weak lower semi-continuity, $u\in L^\infty_{\omega,t,x}(L^1_v)$. Since $\rho_{\varepsilon_1} \in  {\mathcal C}_t(L^1_{\omega,x})$, let us check that: $\rho\in{\mathcal C}_t(L^1_{\omega,x})$. To show this result, one first supposes that $\rho_0\in BV_x$. Since for almost all $\omega$, every $v\in\mathbb{R}$, $X^v(s,t,\cdot)$ is a homeomorphism from $\mathbb{R}^d$ to $\mathbb{R}^d$, (\ref{3.10}) and (\ref{3.11}) are equivalent. If one replaces $f^\prime(v)$ by $f^\prime(\tilde{v})$ in (\ref{3.1}), in view of (\ref{3.10}), one gains that
\begin{eqnarray*}
u_{\varepsilon_1}(t,X^{\tilde{v}}(0,t,x),v)=\frac{1}{\varepsilon_1}\int\limits^t_0e^{\frac{s-t}{\varepsilon_1}}
\chi_{\rho_{\varepsilon_1}(s,X^{\tilde{v}}(0,s,x))}(v)ds+
e^{-\frac{t}{\varepsilon_1}}\chi_{\rho_0(x)}(v),
\end{eqnarray*}
which suggests that
\begin{eqnarray*}
\rho_{\varepsilon_1}(t,X^{\tilde{v}}(0,t,x))=\frac{1}{\varepsilon_1}
\int\limits^t_0e^{\frac{s-t}{\varepsilon_1}}
\rho_{\varepsilon_1}(s,X^{\tilde{v}}(0,s,x))ds+
e^{-\frac{t}{\varepsilon_1}}\rho_0(x).
\end{eqnarray*}
Therefore, we have
\begin{eqnarray}\label{3.25}
\rho_{\varepsilon_1}(t,X^{\tilde{v}}(0,t,x))=\rho_0(x).
\end{eqnarray}
If one set $\tilde{\rho}_{\varepsilon_1}^{\tilde{v}}(t,x)=\rho_{\varepsilon_1}(t,X^{\tilde{v}}(0,t,x))$, then $\tilde{\rho}_{\varepsilon_1}^{\tilde{v}}(\cdot,\cdot)\in L^\infty_{\omega,t}(BV_x)$.

From (\ref{3.25}), for every $R>0$, then
\begin{eqnarray}\label{3.26}
&&\int\limits_{B_R}|\nabla_x\rho_{\varepsilon_1}(t,x)|dx\cr\cr&=&
\int\limits_{B_R}|\nabla_x(\tilde{\rho}_{\varepsilon_1}^{\tilde{v}}
(t,X^{\tilde{v}}(t,0,x)))|dx
\cr\cr&=&\int\limits_{B_R}|\nabla_{X^{\tilde{v}}}
\tilde{\rho}_{\varepsilon_1}^{\tilde{v}}
(t,X^{\tilde{v}}(t,0,x))\nabla_xX^{\tilde{v}}(t,0,x)|dx
\cr\cr&\leq&\int\limits_{B_R}|\nabla_{X^{\tilde{v}}}
\tilde{\rho}_{\varepsilon_1}^{\tilde{v}}
(t,X^{\tilde{v}}(t,0,x))|dx\sup_{(t,x)\in [0,T]\times B_R}|\nabla_xX^{\tilde{v}}(t,0,x)|
\cr\cr&\leq&\int\limits_{\mathbb{R}^d}|\nabla_x\tilde{\rho}_{\varepsilon_1}^{\tilde{v}}
(t,x)|\exp(-f^\prime(v)\int\limits^t_0\mbox{div} b(X^{\tilde{v}}(r,0,x))dr)
\sup_{(t,x)\in [0,T]\times B_R}|\nabla_xX^{\tilde{v}}(t,0,x)|\cr\cr&\leq&  e^{C_0T}\|\rho_0\|_{BV_x}\sup_{(t,x)\in [0,T]\times B_R}|\nabla_xX^{\tilde{v}}(t,0,x)|.
\end{eqnarray}
where in the fifth line we have used the Euler's identity, and the constant $C_0$ is given in (\ref{3.19}). Hence, from the discussion in Step 1, one concludes that
$\rho_{\varepsilon_1}\in L^p(\Omega;L^\infty((0,T);BV_x(B_R)))$ for every $p,R\geq 1$ and particularly, for almost all $\omega\in\Omega$, $\rho_{\varepsilon_1}\in L^\infty((0,T);BV_{loc})$.

Let $\{e_j, j=1,2,_{\cdot},d\}$ denote the standard orthogonal basis of $\mathbb{R}^d$. If one substitutes $f^\prime(\tilde{v})$ for $f^\prime(v)$ in (\ref{3.1}), then for every $1\leq j\leq d$, in view of (\ref{3.8}),
\begin{eqnarray*}
u_{\varepsilon_1}(t,X^{\tilde{v}}(0,t,x+he_j),v)&=&\frac{1}{\varepsilon_1}\int\limits^t_0e^{\frac{s-t}{\varepsilon_1}}
\chi_{\rho_{\varepsilon_1}(s,X^{\tilde{v}}(0,s,x+he_j))}(v)ds+
e^{-\frac{t}{\varepsilon_1}}\chi_{\rho_0(x+he_j)}(v)\cr\cr&=&
\frac{1}{\varepsilon_1}\int\limits^t_0e^{\frac{s-t}{\varepsilon_1}}
\chi_{\tilde{\rho}_{\varepsilon_1}^{\tilde{v}}(s,x+he_j)}(v)ds+
e^{-\frac{t}{\varepsilon_1}}\chi_{\rho_0(x+he_j)}(v).
\end{eqnarray*}
By (\ref{3.25}), therefore
\begin{eqnarray}\label{3.27}
&&\|u_{\varepsilon_1}(t,X^{\tilde{v}}(0,t,\cdot+he_j),\cdot)-
u_{\varepsilon_1}(t,X^{\tilde{v}}(0,t,\cdot,\cdot)\|_{L^1_{x,v}}
\cr\cr&\leq&
\int\limits_{\mathbb{R}^d_x\times\mathbb{R}_v}  \frac{1}{\varepsilon_1}\int\limits^t_0e^{\frac{s-t}
{\varepsilon_1}}|
\chi_{\tilde{\rho}_{\varepsilon_1}^{\tilde{v}}(s,x+he_j)}(v)-
\chi_{\tilde{\rho}_{\varepsilon_1}^{\tilde{v}}(s,x)}(v)|dxdvds \cr\cr &&+
\int\limits_{\mathbb{R}^d_x\times\mathbb{R}_v}e^{-\frac{t}{\varepsilon_1}}|\chi_{\rho_0(x+he_j)}(v)-
\chi_{\rho_0(x)}(v)|dxdv
\cr\cr&\leq&
\int\limits_{\mathbb{R}^d_x}  \frac{1}{\varepsilon_1}\int\limits^t_0e^{\frac{s-t}
{\varepsilon_1}}|\tilde{\rho}_{\varepsilon_1}^{\tilde{v}}(s,x+he_j)-
\tilde{\rho}_{\varepsilon_1}^{\tilde{v}}(s,x)|dxds +
\int\limits_{\mathbb{R}^d_x}e^{-\frac{t}{\varepsilon_1}}|\rho_0(x+he_j)-
\rho_0(x)|dx
\cr\cr &\leq& \Big[
  \frac{1}{\varepsilon_1}\int\limits^t_0e^{\frac{s-t}
{\varepsilon_1}}ds+
e^{-\frac{t}{\varepsilon_1}}\Big]|h|\|\rho_0\|_{BV_x}
\cr\cr
&=&|h|\|\rho_0\|_{BV_x}, \quad {\mathbb P}-a.s.,
\end{eqnarray}

With the help of Theorem 1.7.2 (see \cite[P17]{Daf}), from (\ref{3.27}),
then $u_{\varepsilon_1}(t,X^{\tilde{v}}(0,t,\cdot),\cdot)\in BV_x(L^1_v)$. Let $C_0$ be given in (\ref{3.19}), by an analogue manipulation of (\ref{3.26}), we derive that:
\begin{eqnarray}\label{3.28}
\int\limits_{B_R}\int\limits_{\mathbb{R}_v}|\nabla_x
u_{\varepsilon_1}(t,x,v)|dvdx
\leq  e^{C_0T}\|\rho_0\|_{BV_x}\sup_{(t,x)\in [0,T]\times B_R}|\nabla_xX^{\tilde{v}}(t,0,x)|.
\end{eqnarray}
So $u_{\varepsilon_1}\in L^p(\Omega;L^\infty((0,T);BV_x(B_R;L^1(\mathbb{R}_v))))$ for every $p,R\geq 1$ and particularly, for almost all $\omega\in\Omega$, $u_{\varepsilon_1}\in L^\infty((0,T);BV_{loc}(\mathbb{R}^d_x;L^1(\mathbb{R}_v)))$.

On the other hand, if we define $u^1_{\varepsilon_1}(t,x,v)=u_{\varepsilon_1}(t,x+B(t),v)$, then (\ref{3.1}) has an equivalent representation below (the proof is similar to (\ref{3.5})-(\ref{3.12})):
\begin{eqnarray}\label{3.29}
\left\{
  \begin{array}{ll}
\partial_tu^1_{\varepsilon_1}(t,x,v)+
f^\prime(v)b_1(x)\cdot\nabla_x u^1_{\varepsilon_1}
= \frac{1}{\varepsilon_1}
{\big[} \chi_{\rho^1_{\varepsilon_1}}-u^1_{\varepsilon_1}{\big]}, \\
u^1_{\varepsilon_1}(t,x,v)|_{t=0}=\chi_{\rho_0(x)}(v), \\ \rho^1_\varepsilon(t,x)=\int_{\mathbb{R}} u^1_{\varepsilon_1}(t,x,v)dv, \\ b_1(x)=b(x+B(t)).
  \end{array}
\right.
\end{eqnarray}

From (\ref{3.29}), with the help of inequality (\ref{3.28}), for every $0\leq t_1\leq t_2\leq T$, every $R>0$, we fulfill
\begin{eqnarray}\label{3.30}
&&\|\rho_{\varepsilon_1}(t_1,\cdot+B(t_1))-\rho_{\varepsilon_1}
(t_2,\cdot+B(t_2))\|_{L^1(B_R)}
\cr\cr&=&\int\limits_{t_1}^{t_2}{\Big \|}\sum_{i=1}^d b_2(\cdot)\cdot\nabla_x\int\limits_{\mathbb{R}}f^\prime(v)u^1_{\varepsilon_1}(t,\cdot,v)dv
{\Big \|}_{L^1(B_R)}dt
\cr\cr&\leq& Ce^{C_0T}\|b\|_{L^\infty_x}\|f^\prime\|_{L^\infty_v}|t_1-t_2|
\|\rho_0\|_{BV_x}, \quad {\mathbb P}-a.s..
\end{eqnarray}
By (\ref{3.26}) and (\ref{3.30}), so
\begin{eqnarray}\label{3.31}
&&\|\rho_{\varepsilon_1}(t_1)-\rho_{\varepsilon_1}(t_2)\|_{L^1(B_R)}
\cr\cr&\leq&
\|\rho_{\varepsilon_1}(t_2,\cdot)-\rho_{\varepsilon_1}
(t_2,\cdot+B(t_2)-B(t_1))\|_{L^1(B_R)}
+ C|t_1-t_2|\|\rho_0\|_{BV_x}
\cr\cr&\leq&C(|B(t_2)-B(t_1)| +|t_1-t_2|)\|\rho_0\|_{BV_x}, \quad {\mathbb P}-a.s..
\end{eqnarray}
By virtue of (\ref{2.19}) in Lemma \ref{lem2.2},  from (\ref{3.31}), we arrive at
\begin{eqnarray}\label{3.32}
&&\|\rho_{\varepsilon_1}(t_1,\cdot)-\rho_{\varepsilon_1}(t_2,\cdot)\|_{L^1(B_R)}
\cr\cr&\leq& C(\sqrt{|t_1-t_2|\log\log(1/|t_1-t_2|)} +|t_1-t_2|)\|\rho_0\|_{BV_x}
\cr\cr&\leq& C(|t_1-t_2|^{1/3}+|t_1-t_2|)\|\rho_0\|_{BV_x}, \quad  {\mathbb P}-a.s.,
\end{eqnarray}
if $|t_1-t_2|$ is sufficiently small.

Combining (\ref{3.26}) and (\ref{3.32}), with the aid of Helly's theorem (see \cite[P17]{Daf}) and Ascoli-Arzela's compact criterion, given $T^0<T$ (sufficiently small), then for almost all $\omega\in\Omega$,
\begin{eqnarray*}
\rho_{\varepsilon_1}(\cdot,\cdot,\omega)\longrightarrow \rho^\omega\in {\mathcal C}([0,T^0];L^1_{loc}(\mathbb{R}^d)), \  \mbox{as} \ \varepsilon_1\rightarrow 0.
\end{eqnarray*}
After repeating above calculations finitely many times, one derives that for almost all $\omega\in\Omega$
\begin{eqnarray}\label{3.33}
\rho_{\varepsilon_1}(\cdot,\cdot,\omega)\longrightarrow \rho^\omega\in {\mathcal C}([0,T];L^1_{loc}(\mathbb{R}^d)), \  \mbox{as} \ \varepsilon_1\rightarrow 0.
\end{eqnarray}
According to (\ref{3.20}) and (\ref{3.33}), by applying dominated convergence, one concludes
\begin{eqnarray}\label{3.34}
\rho_{\varepsilon_1}(\cdot,\cdot,\cdot)\longrightarrow \rho(\cdot,\cdot,\cdot)\in L^1(\Omega;{\mathcal C}([0,T];L^1_{loc}(\mathbb{R}^d))), \  \mbox{as} \ \varepsilon_1\rightarrow 0.
\end{eqnarray}

On the other hand, by applying Fatou's lemma and inequality (\ref{3.19}), for every $R>0$, it yields that
\begin{eqnarray}\label{3.35}
{\mathbb E}\sup_{0\leq t\leq T}\int\limits_{B_R}|\rho^\omega(t,x)|dx&\leq& {\mathbb E}\liminf_{\varepsilon_1\rightarrow0}\sup_{0\leq t\leq T}\int\limits_ {B_R}|\rho^\omega_{\varepsilon_1}(t,x)|dx\cr\cr
&\leq&\liminf_{\varepsilon_1\rightarrow0}\sup_{0\leq t\leq T}\int\limits_{\mathbb{R}^{d+1}_{x,v}}|u_{\varepsilon_1}(t,x,v,\omega)|dxdv
\cr\cr&\leq&e^{C_0T}\|\rho_0\|_{L^1_x},\ \ \ {\mathbb P}-a.s..
\end{eqnarray}
So $\rho\in{\mathcal C}_t(L^1_{\omega,x})$ and this completes the proof for $BV$ initial data. For general initial value $\rho_0\in L^\infty\cap L^1(\mathbb{R}^d)$, we can justify it by standard cutoff and $BV$-regularization of initial data (consult to \cite{CM}), we omit the details here.

 \vskip2mm\noindent
\textbf{Step 5.} $\frac{1}{\varepsilon_1} {\big[}
\chi_{\rho_{\varepsilon_1}}-u_{\varepsilon_1}{\big]}=\partial_v
m_{\varepsilon_1}$, $u(t,x,v)=\chi_{\rho(t,x)}(v)$ and
$\rho$ solves (\ref{1.6}).

Let $(t,x)\in (0,T)\times \mathbb{R}^d$ be fixed, assuming without
loss of generality that $\rho_{\varepsilon_1}\geq0$, define
\begin{eqnarray*}
 m_{\varepsilon_1}(t,x,v)= \frac{1}{\varepsilon_1}
\int\limits_{-\infty}^v[
\chi_{\rho_{\varepsilon_1}(t,x)}(r)-u_{\varepsilon_1}(t,x,r)]dr.
\end{eqnarray*}
In view of (\ref{3.11}),
$$
u_{\varepsilon_1}(t,x,r)\in \cases{ \ \ [0,1], \ \ \mbox{when} \ r>0,
\cr  \ [-1,0], \ \mbox{when} \ r<0.}
$$
Hence $m_{\varepsilon_1}(t,x,v)$ is nondecreasing on
$(-\infty,\rho_{\varepsilon_1})$ and nonincreasing on
$[\rho_{\varepsilon_1},\infty)$. Since $m_{\varepsilon_1}(t,x,-\infty)=m_\varepsilon(t,x,\infty)=0$, we
conclude $m_{\varepsilon_1}\geq 0$.

Observing that $\rho_0\in L^\infty\cap L^1(\mathbb{R}^d)$, owing to (\ref{3.11}) and (\ref{3.20}), we get
$$
\mbox{supp}m_{\varepsilon_1}\subset [0,T]\times \mathbb{R}^d\times [-N,N],
$$
where $N=\|\rho_0\|_{L^\infty(\mathbb{R}^d)}$.

Repeating the calculations from (\ref{2.9}) to (\ref{2.14}), it yields that
\begin{eqnarray}\label{3.36}
&& {\mathbb E}{\Big |}\int\limits^T_0dt\int\limits_{\mathbb{R}^d}dx\int\limits_{\mathbb{R}}
m_{\varepsilon_1}(t,x,v)dv{\Big |}^2\cr\cr&=&
{\mathbb E}{\Big |}\int\limits^T_0dt\int\limits_{\mathbb{R}^d}dx\int\limits_{-N}^Ndv\int\limits_{-N}^v
[-\partial_tu_{\varepsilon_1}(t,x,r)+\mbox{div} b(x)f^\prime(r)u_{\varepsilon_1}(t,x,r)]dr{\Big |}^2 \cr \cr &\leq&
12N^2\Big[{\mathbb E}\|u_{\varepsilon_1}(T)\|^2_{L^1_{x,v}}+
\|\rho_0\|_{L^1_x}^2+\|\mbox{div} b\|^2_{L^\infty_x}\|f^\prime\|^2_{L^\infty_v}{\mathbb E}{\Big |}\int\limits_0^T
\|u_{\varepsilon_1}(t)\|_{L^1_{x,v}}dt{\Big |}^2\Big].
\end{eqnarray}
Thanks to (\ref{3.19}), it follows from (\ref{3.36}) that
\begin{eqnarray*}
{\mathbb E}{\Big |}\int\limits^T_0dt\int\limits_{\mathbb{R}^d}dx\int\limits_{\mathbb{R}}
m_{\varepsilon_1}(t,x,v)dv{\Big |}^2
\leq 12N^2\Big[e^{2C_0T}+1+C^2_0T^2e^{2C_0T}\Big]
\|\rho_0\|^2_{L^1_x}.
\end{eqnarray*}
Whence $m_{\varepsilon_1}$ is bounded uniformly in $\varepsilon_1$.

By extracting a unlabeled subsequence, one achieves
$$
m_{\varepsilon_1}\stackrel{w*}{\longrightarrow} m \geq 0 \ \mbox{in} \ \ L^2(\Omega;{\mathcal M}_b([0,T]\times\mathbb{R}^{d+1}_{x,v})), \ \mbox{as} \ \varepsilon_1\rightarrow 0.
$$
In particular, for every $\varphi\in {\mathcal D}(\mathbb{R}^d)$, every $\psi\in {\mathcal D}(\mathbb{R})$, every $h\in L^2(\Omega)$ and every $l\in L^2(0,T)$,
$$
\lim_{\varepsilon_1\rightarrow0}{\mathbb E}\Big[h\int\limits_0^Tl(t)m_{\varepsilon_1}
(\varphi\psi)(t)dt\Big]
={\mathbb E}\Big[h\int\limits_0^Tl(t)m(\varphi\psi)(t)dt\Big]
$$
where
$$
m_{\varepsilon_1}(\varphi\psi)(t)=\int\limits_0^t\int\limits_{\mathbb{R}^{d+1}_{x,v}}\varphi(x)\psi(v)
m_{\varepsilon_1}(ds,dx,dv),  $$
hence (\ref{2.3}) holds.

Observing that $u_{\varepsilon_1}\stackrel{w*}{\longrightarrow} u \ \mbox{in} \ L^\infty_{\omega,t,x,v}$, $\rho_{\varepsilon_1}\longrightarrow \rho$ in ${\mathcal C}([0,T];L^1(\Omega;L^1_{loc}(\mathbb{R}^d)))$
and $m_{\varepsilon_1}\stackrel{w*}{\longrightarrow} m$ in $L^2(\Omega;{\mathcal M}_b([0,T]\times\mathbb{R}^{d+1}_{x,v}))$,
so $u_{\varepsilon_1}(t,x,v)-\chi_{\rho_{\varepsilon_1}}(v)\rightarrow 0$ in the distributional sense. Then $u=\chi_{\rho(t,x)}(v)\in{\mathcal C}_t(L^1_{\omega,x,v})$  meets (\ref{2.3}) in the distributional sense, for almost all $\omega\in\Omega$. On the other hand, $u\in{\mathcal C}_t(L^1_{\omega,x,v})$, which implies (\ref{2.5}) holds, and by Lemma \ref{lem2.1}, (\ref{2.4}) is true, thus $\rho$ is a stochastic entropy solution to (\ref{1.6}). $\Box$

\vskip2mm\noindent
\textbf{Remark 3.1.} (i) Consider the multi-dimensional scalar conservation laws driven by a multiplicative cylindrical Brownian motion $W(t)$ on a $d$-dimensional torus
\begin{eqnarray}\label{3.37}
\partial_t\rho(t,x)+\mbox{div}_xF(\rho)
=A(\rho)\dot{W}(t), \quad  t\in(0,T), \ \ x\in {\mathbb T}^d.
\end{eqnarray}
Under the assumption that: $F\in {\mathcal C}^{4,\eta}$ for some $\eta>0$,  with a polynomial growth of its first derivative, $\rho_0\in L^p(\Omega\times {\mathbb T}^d)$ for all $p\in [1,\infty)$, in \cite{Hof},
using stochastic BGK approximation, Hofmanov\'{a} obtained the existence of the stochastic kinetic solution, which lied in $L^p(\Omega;L^\infty([0,T]; L^p({\mathbb T}^d)))$ for all $p\in [1,\infty)$ (for more details about (\ref{3.37}), one also consults to \cite{DV,VW}). Different from \cite{Hof}, under the assumption that $\rho_0\in L^\infty\cap L^1$, we prove that (\ref{1.6}) exists a stochastic entropy solution (it is also a stochastic kinetic solution). Besides the stochastic entropy solution is a continuous semi-martingale.

(ii) For simplicity, here we only discuss the noise given by $\partial_{x_i}\rho\circ\dot{B}_i(t)$, however, the method is appropriate for the stochastic balance law below
\begin{eqnarray*}
\partial_t\rho(t,x)+b(x)\cdot\nabla_x f(\rho)
+\partial_{x_i}G_{i,j}(t,\rho)\circ\dot{B}_j(t)=A(t,\rho), \ (\omega,t,x)\in \Omega\times(0,T)\times\mathbb{R}^d,
\end{eqnarray*}
here $1\leq i\leq  d, \ 1\leq j\leq n$, $d,n\in{\mathbb N}$. But now one should replace the SDE (\ref{3.3}) by
\begin{eqnarray*}
\left(
\begin{array}{cc}
dX(t) \\ dV(t)
  \end{array}
\right)=\left(
\begin{array}{cc}
 f^\prime(v)b(X(t))dt+\partial_vG(t,V)dB^\top(t)  \\
 A(t,V)dt
\end{array} \right), \ \
\left(\begin{array}{cc}
 X(0) \\ V(0)
  \end{array}
\right)=\left(
\begin{array}{cc}
 x  \\ v
  \end{array}
\right),
 \end{eqnarray*}
where $G(t,v)=(G_{i,j}(t,v))\in \mathbb{R}^{d\times n}$, $B=(B_1,B_2,_{\cdots}, B_n)$ is a standard $n$-dimensional Brownian motion.

\section{Uniqueness}\label{sec4}
\setcounter{equation}{0}
This section is devoted to prove the second part of Theorem \ref{the2.1}. That is stated as below.

\begin{theorem} \label{the4.1} Let $\rho_0\in L^\infty\cap L^1(\mathbb{R}^d)$, $b\in {\mathcal C}_b^\alpha(\mathbb{R}^d;\mathbb{R}^d)\cap BV_{loc}(\mathbb{R}^d;\mathbb{R}^d)$ with $\alpha\in (0,1)$, $f\in {\mathcal C}^1(\mathbb{R})$. We suppose in addition that $\mbox{div}_xb\in L^\infty(\mathbb{R}^d)$. Then the stochastic entropy solution to the Cauchy problem (\ref{1.6}) is unique.
 \end{theorem}
\vskip1mm\noindent
\textbf{Proof.} Let $\rho_1$ and $\rho_2$ be two stochastic entropy solutions of (\ref{1.6}), with initial values $\rho_{0,1}$ and
$\rho_{0,2}$, respectively. Then $u_1=\chi_{\rho_1}$ and
$u_2=\chi_{\rho_2}$ are stochastic weak solutions
of (\ref{2.2}) with nonhomogeneous terms $\partial_vm_1$ and
$\partial_vm_2$, initial data $u_{0,1}=\chi_{\rho_{0,1}}$ and
$u_{0,2}=\chi_{\rho_{0,2}}$, respectively.

Suppose that $\varrho_1$ is a regularization kernel stated in Lemma \ref{lem2.3}. Let $\varrho$ be another regularization kernel, i.e. $\varrho \in {\mathcal D}_{+}(\mathbb{R}),  \ \int_{\mathbb{R}} \varrho(r)dr=1$. For $\varepsilon_2>0$, set $\varrho_{\varepsilon_2}(t)=\frac{1}{\varepsilon_2}\varrho
(\frac{t}{\varepsilon_2})$, then
$u_\iota^{\varepsilon}:=u_\iota\ast\varrho_{1,\varepsilon_1}\ast
\varrho_{\varepsilon_2}$ ($\iota=1,2$)
meet
\begin{eqnarray*}
\left\{
  \begin{array}{ll}
\partial_tu_\iota^{\varepsilon}+f^\prime(v)b(x)\cdot \nabla_{x}
u_\iota^{\varepsilon}+\partial_{x_i}u_\iota^{\varepsilon}
\circ\dot{B}_i(t)
=\partial_vm_\iota^{\varepsilon}+
R_\iota^{\varepsilon},  \\
u_\iota^{\varepsilon}(t,x,v)|_{t=0}=\chi_{\rho_0^\iota}
\ast\varrho_{1,\varepsilon_1}(x,v),
  \end{array}
\right.
\end{eqnarray*}
where
\begin{eqnarray*}
u_\iota(t,x,v):=0, \ \mbox{when} \ t\bar{\in}[0,T], \   R_{\iota}^{\varepsilon}=f^\prime(v)b(x)\cdot \nabla_{x}
u_\iota^{\varepsilon}-f^\prime(v)[b(x)\cdot \nabla_{x}
u_\iota]^{\varepsilon}.
\end{eqnarray*}

Define
\begin{eqnarray}\label{4.1}
\sigma_{\varepsilon}:&=&
\partial_t[|u_1^{\varepsilon}|+|u_2^{\varepsilon}|
-2u_1^{\varepsilon}u_2^{\varepsilon}]
+f^\prime(v)b(x)\cdot \nabla_{x}[|u_1^{\varepsilon}|+|u_2^{\varepsilon}|
-2u_1^{\varepsilon}u_2^{\varepsilon}]
\cr\cr&&+
\partial_{x_i}[|u_1^{\varepsilon}|+|u_2^{\varepsilon}|
-2u_1^{\varepsilon}u_2^{\varepsilon}]
\circ\dot{B}_i
-\mbox{sgn}(v)\partial_v(m_1^{\varepsilon}+m_2^{\varepsilon})
\cr\cr&&+2(u_1^{\varepsilon}\partial_vm_2^{\varepsilon}+
u_2^{\varepsilon}\partial_vm_1^{\varepsilon}),
\end{eqnarray}
then
\begin{eqnarray}\label{4.2}
\sigma_{\varepsilon}=\mbox{sgn}(v)(R_1^{\varepsilon}+R_2^\varepsilon)
-2(u_1^{\varepsilon}R_2^{\varepsilon}+u_2^{\varepsilon}R_1^{\varepsilon}),
\end{eqnarray}
for $|u_\iota^{\varepsilon}|=\mbox{sgn}(v)u_\iota^{\varepsilon}$.

Since $u_1^{\varepsilon}$ and $u_2^{\varepsilon}$ are bounded in $x,v$, using the estimates (\ref{2.20}), (\ref{2.21}) in Lemma \ref{lem2.3}, then for every $K=K_1\times K_2$ with $K_1\subset \subset \mathbb{R}^d, K_2\subset \subset \mathbb{R}$, we derive
\begin{eqnarray}\label{4.3}
&&\limsup_{\varepsilon_1\downarrow0}\limsup_{\varepsilon_2\downarrow0}\int\limits_K |R_\iota^{\varepsilon}|dxdv
\cr\cr&\leq& {\mathcal L}(K_2) \int\limits_{K_1}
\Lambda(M(x),\varrho_1(x))d|D^sb|(x)+C{\mathcal L}(K_2)(d + I(\varrho_1))|D^ab|(K_1),
\end{eqnarray}
and
\begin{eqnarray}\label{4.4}
\limsup_{\varepsilon_1\downarrow0}\limsup_{\varepsilon_2\downarrow0}
\int\limits_K |R_\iota^{\varepsilon}|dxdv \leq {\mathcal L}(K_2)I(\varrho_1)|D^s
b|(K_1),
\end{eqnarray}
for almost all $(\omega,t)\in \Omega\times [0,T]$, with $\iota=1,2$.

 Thanks to (\ref{4.2})-(\ref{4.4}), for $\varepsilon_1$ sufficiently small, let $K$ be given above, then $\sigma_{\varepsilon}\in L^\infty(\Omega\times[0,T];{\mathcal M}_b(K))$, and there is a constant $C>0$, such that
\begin{eqnarray}\label{4.5}
&&\limsup_{\varepsilon_1\downarrow0}\limsup_{\varepsilon_2\downarrow0}\int\limits_K |\sigma_{\varepsilon}|dxdv
\cr\cr&\leq& C{\mathcal L}(K_2) \int\limits_{K_1}
\Lambda(M(x),\varrho_1(x))d|D^sb|(x)+C{\mathcal L}(K_2)(d + I(\varrho_1))|D^ab|(K_1),
\end{eqnarray}
and
\begin{eqnarray}\label{4.6}
\limsup_{\varepsilon_1\downarrow0}\limsup_{\varepsilon_2\downarrow0}
\int\limits_K |\sigma_{\varepsilon}|dxdv \leq C{\mathcal L}(K_2)I(\varrho_1)|D^s
b|(K_1),
\end{eqnarray}
for almost all $(\omega,t)\in \Omega\times [0,T]$, which hint that
$$
\sigma_{\varepsilon}\stackrel{w*}{\longrightarrow} \sigma \ \mbox{in} \ \ L^\infty(\Omega\times[0,T];{\mathcal M}_b(K)), \ \mbox{as} \ \varepsilon_2\rightarrow 0 \ \mbox{first}, \ \varepsilon_1\rightarrow 0 \ \mbox{next}.
$$

Now we show $\sigma=0$. Indeed, with the aid of (\ref{4.6}), one concludes $\sigma$ is singular with Lebesgue measure ${\mathcal L}^{d+1}$, then (\ref{4.5}) uses, we derive for every $\varphi \in {\mathcal D}(\mathbb{R}^d), \psi\in{\mathcal D}(\mathbb{R})$,
\begin{eqnarray}\label{4.7}
|\langle\sigma,\varphi\psi\rangle| \leq C \int\limits_{\mathbb{R}^d}
\Lambda(M(x),\varrho_1(x))d|D^sb|(x)\int\limits_{\mathbb{R}}|\psi(v)|dv, \ a.s. \ (\omega,t)\in\Omega\times[0,T].
\end{eqnarray}
From (\ref{4.7}) and Lemma \ref{lem2.3} (ii), by minimizing the even kernel, we claim that $\sigma=0$.

Let $\theta$ and $\theta_n$ be given in (\ref{2.11}), $\xi$ and $\xi_k$ be given in (\ref{2.16}). From (\ref{4.2}), for almost all $(\omega,t)\in\Omega\times[0,T]$, it follows that
\begin{eqnarray}\label{4.8}
&&\int\limits_{\mathbb{R}^{d+1}_{x,v}}\sigma_{\varepsilon}\xi_k(v)
\theta_n(x)dxdv
\cr\cr&=&\frac{d}{dt} \int\limits_{\mathbb{R}^{d+1}_{x,v}}[|u_1^{\varepsilon}|+|u_2^{\varepsilon}|
-2u_1^{\varepsilon}u_2^{\varepsilon}]\xi_k(v)
\theta_n(x)dxdv\cr\cr&&-
\int\limits_{\mathbb{R}^{d+1}_{x,v}}[|u_1^{\varepsilon}|+|u_2^{\varepsilon}|
-2u_1^{\varepsilon}u_2^{\varepsilon}]
\xi_k(v)f^\prime(v)b(x)\cdot\nabla_x\theta_n(x)dxdv\cr\cr&&-
\int\limits_{\mathbb{R}^{d+1}_{x,v}}[|u_1^{\varepsilon}|+|u_2^{\varepsilon}|
-2u_1^{\varepsilon}u_2^{\varepsilon}]
\partial_{x_i}\theta_n(x)\xi_k(v)\circ\dot{B}_i(t)dxdv +J,
\end{eqnarray}
where
\begin{eqnarray*}
J&=&\int\limits_{\mathbb{R}^{d+1}_{x,v}}\xi_k\theta_n\mbox{sgn}(v)\partial_v(m_1^{\varepsilon}+
m_2^{\varepsilon})dxdv-2\int\limits_{\mathbb{R}^{d+1}_{x,v}}\xi_k\theta_n[
u_1^{\varepsilon}\partial_vm_2^{\varepsilon}+
u_2^{\varepsilon}\partial_vm_1^{\varepsilon}]dxdv
\cr\cr&=:&J_1-2J_2.
\end{eqnarray*}

Obviously, $m_\iota ^{\varepsilon} \ (\iota=1,2)$ is continuous in $v$ in a neighborhood of zero, therefore for large $k$,
\begin{eqnarray}\label{4.9}
J_1=-2\int\limits_{\mathbb{R}^d}\theta_n(x)[
m_1^{\varepsilon}(t,x,0)+m_2^{\varepsilon}(t,x,0)]dx.
\end{eqnarray}

Moreover, noting the fact $m_\iota\geq 0 \ (\iota=1,2)$,
so for $k$ large enough,
\begin{eqnarray}\label{4.10}
J_2&=&\int\limits_{\mathbb{R}^{d+1}_{x,v}}\xi_k\theta_n[
u_1^{\varepsilon}\partial_vm_2^{\varepsilon}+
u_2^{\varepsilon}\partial_vm_1^{\varepsilon}]dxdv
\cr\cr
&=&\int\limits_{\mathbb{R}^d}\theta_n(x)dx\int\limits_{\mathbb{R}}\Big[
\int\limits_{\mathbb{R}^{d+1}}\chi_{\rho_1(s,y)}(v)\varrho_{1,\varepsilon_1}(x-y)
\varrho_{\varepsilon_2}(t-s)dsdy\partial_vm_2^{\varepsilon}(t,x,v)\cr\cr&&+
\int\limits_{\mathbb{R}^{d+1}}\chi_{\rho_2(s,y)}(v)\varrho_{1,\varepsilon_1}(x-y)
\varrho_{\varepsilon_2}(t-s)dsdy\partial_vm_1^{\varepsilon}(t,x,v) \Big]dv
\cr\cr
&=&\int\limits_{\mathbb{R}^d}\theta_n(x)dx\int\limits_{\mathbb{R}}{\Big\{}
\int\limits_{\mathbb{R}^{d+1}}
\chi_{\rho_1(s,y)}(v) [\delta_{\rho_1(s,y)}(v)-\delta_{0}(v)] \varrho_{1,\varepsilon_1}(x-y)
\varrho_{\varepsilon_2}(t-s)dsdym_2^{\varepsilon}(t,x,v)\cr\cr&&+
\int\limits_{\mathbb{R}^{d+1}}[\delta_{\rho_2(s,y)}(v)-\delta_{0}(v)] \varrho_{1,\varepsilon_1}(x-y)
\varrho_{\varepsilon_2}(t-s)dsdy m_1^{\varepsilon}(t,x,v) {\Big\}}dv
\cr\cr&\geq& -\int\limits_{\mathbb{R}^d}\theta_n(x)[m_1^{\varepsilon}(t,x,0)+
m_2^{\varepsilon}(t,x,0)]dx.
\end{eqnarray}

Combining (\ref{4.9})-(\ref{4.10}), for $n$ and $k$ ($k$ is big enough) be fixed, if one lets $\varepsilon_2$ tend to zero first,
$\varepsilon_1$ incline to zero next in (\ref{4.8}), it yields \begin{eqnarray}\label{4.11}
&&\int\limits_{\mathbb{R}^{d+1}_{x,v}}[|u_1|+|u_2|
-2u_1u_2]\xi_k(v) \theta_n(x)dxdv
 \cr\cr &\leq&
\int\limits_0^t\int\limits_{\mathbb{R}^{d+1}_{x,v}}[|u_1|+|u_2|
-2u_1u_2]
\xi_kf^\prime(v)\mbox{div}_x(b(x)\theta_n(x))dxdvds\cr\cr&&+
\int\limits_0^t\circ dB_i(s)\int\limits_{\mathbb{R}^{d+1}_{x,v}}[|u_1|+|u_2|
-2u_1u_2]
\partial_{x_i}\theta_n\xi_kdxdv
\cr\cr &=&
\int\limits_0^t\int\limits_{\mathbb{R}^{d+1}_{x,v}}[|u_1|+|u_2|
-2u_1u_2]
\xi_kf^\prime(v)\mbox{div}_x(b(x)\theta_n(x))dxdvds\cr\cr&&
-\frac{1}{2}
\int\limits_0^tds\int\limits_{\mathbb{R}^{d+1}_{x,v}}[|u_1|+|u_2|
-2u_1u_2]
\Delta_x\theta_n\xi_kdxdv\cr\cr&&+\int\limits_0^tdB_i(s)\int\limits_{\mathbb{R}^{d+1}_{x,v}}[|u_1|+|u_2|
-2u_1u_2]
\partial_{x_i}\theta_n\xi_kdxdv,\quad {\mathbb P}-a.s.,
\end{eqnarray}
where in the last identity, we used partial integration and the following fact (see \cite{FGP}),
\begin{eqnarray*}
&&\partial_{x_i}[|u_1|+|u_2|
-2u_1u_2]\circ dB_i(t)
\cr\cr&=&\partial_{x_i}[|u_1|+|u_2|
-2u_1u_2] dB_i(t)-\frac{1}{2}\Delta_x[|u_1|+|u_2|
-2u_1u_2]dt, \quad {\mathbb P}-a.s..
\end{eqnarray*}

Since the last term in (\ref{4.11}) is an $\mathcal{F}_t$-martingale, and $|u_1|=|u_1|^2, |u_2|=|u_2|^2, |u_1-u_2|=|u_1-u_2|^2$, from (\ref{4.11}), one concludes
\begin{eqnarray*}
{\mathbb E}\int\limits_{\mathbb{R}^{d+1}_{x,v}}|u_1(t)-u_2(t)|dxdv &\leq&
{\mathbb E}\int\limits_0^t\int\limits_{\mathbb{R}^{d+1}_{x,v}}|u_1(s)-u_2(s)|f^\prime(v)\mbox{div}_xb(x)dxdvds
\cr\cr&\leq&\|f^\prime\|_{L^\infty_v}\|\mbox{div} b\|_{L^\infty_x}
{\mathbb E}\int\limits_0^t\int\limits_{\mathbb{R}^{d+1}_{x,v}}|u_1(s)-u_2(s)|dxdvds,
\end{eqnarray*}
by taking $n$ to infinity first, $k$ to infinity second, which hints
 \begin{eqnarray}\label{4.12}
{\mathbb E}\int\limits_{\mathbb{R}^d}|\rho_1(t)-\rho_2(t)|dx={\mathbb E}\int\limits_{\mathbb{R}^{d+1}_{x,v}}|u_1(t)-u_2(t)|
dxdv\leq0.
\end{eqnarray}
From (\ref{4.12}), we complete the proof. $\Box$

From above calculation, one clearly has the below comparison result.
\begin{corollary}
\textbf{(Comparison Principle)} Let $b,f$ be described in Theorem \ref{the4.1} and $\rho_{0,1},\rho_{0,2}\in L^\infty\cap L^1(\mathbb{R}^d)$.
Assume that $\rho_1$ and $\rho_2$ are two
stochastic entropy solutions of (\ref{1.6}), with initial values
$\rho_{0,1}$ and $\rho_{0,2}$. If $\rho_{0,1}\leq \rho_{0,2}$, then
with probability 1, $\rho_1\leq \rho_2$. In particular, if the initial value is nonnegative, then with probability 1, the unique stochastic solution is nonnegative as well.
\end{corollary}
\textbf{Proof.} Clearly, mimicking above calculation, we have
\begin{eqnarray*}
&&{\mathbb E}\int\limits_{\mathbb{R}^{d+1}_{x,v}}[u_1(t,x,v)-u_2(t,x,v)]dxdv \cr\cr&=&\int\limits_{\mathbb{R}^d}[\rho_{0,1}(x)-\rho_{0,2}(x)]dx+
{\mathbb E}\int\limits_0^t\int\limits_{\mathbb{R}^{d+1}_{x,v}}[u_1(s)-u_2(s)]
f^\prime(v)\mbox{div}_xb(x)dxdvds.
\end{eqnarray*}
Observing that
$$
[u_1(t,x,v)-u_2(t,x,v)]^+=\frac{|u_1-u_2|+(u_1-u_2)}{2},
$$
hence
\begin{eqnarray}\label{4.13}
&&{\mathbb E}\int\limits_{\mathbb{R}^d}[\rho_1(t,x)-\rho_2(t,x)]^+dx
\cr\cr&=&{\mathbb E}\int\limits_{\mathbb{R}^{d+1}_{x,v}}[u_1(t,x,v)-u_2(t,x,v)]^+dxdv
\cr\cr &=&\frac{1}{2}
{\mathbb E}\int\limits_{\mathbb{R}^{d+1}_{x,v}}|u_1(t,x,v)-u_2(t,x,v)|dxdv+\frac{1}{2}
{\mathbb E}\int\limits_{\mathbb{R}^{d+1}_{x,v}}[u_1(t,x,v)-u_2(t,x,v)]dxdv \cr\cr &\leq&
\frac{1}{2}\int\limits_{\mathbb{R}^d}|\rho_{0,1}(x)-\rho_{0,2}(x)|dx+\frac{1}{2}{\mathbb E}\int\limits_0^t\int\limits_{\mathbb{R}^{d+1}_{x,v}}|u_1(s)-u_2(s)|
f^\prime(v)\mbox{div}_xb(x)dxdvds
\cr\cr
&\leq&\int\limits_{\mathbb{R}^d}[\rho_{0,1}(x)-\rho_{0,2}(x)]^+dx+\|f^\prime\|_{L^\infty_v}\|\mbox{div} b\|_{L^\infty_x}
{\mathbb E}\int\limits_0^t\int\limits_{\mathbb{R}^d}[\rho_1(s,x)-\rho_2(s,x)]^+dxds.
\end{eqnarray}
Owing to the Gr\"{o}nwall inequality, one gains from (\ref{4.13}) that
\begin{eqnarray*}
{\mathbb E}\int\limits_{\mathbb{R}^d}[\rho_1(t,x)-\rho_2(t,x)]^+dx
\leq \exp(\|f^\prime\|_{L^\infty_v}\|\mbox{div} b\|_{L^\infty_x}t)\int\limits_{\mathbb{R}^d}[\rho_{0,1}(x)-\rho_{0,2}(x)]^+dx,
\end{eqnarray*}
which implies $\rho_1\leq \rho_2, \ {\mathbb P}-a.s.$. $\Box$

\section{Regularity}\label{sec5}
\setcounter{equation}{0}
In this section, we give the proof for Theorem \ref{the2.2}, some details are described as below.

\vskip2mm\noindent
\textbf{Proof of Theorem \ref{the2.2}.} For $\varepsilon_1>0$ be fixed, if one denotes the unique stochastic weak solution of (\ref{3.1}) by $u_{\varepsilon_1}$, since $\rho_0\in BV_x$, one gets the inequality (\ref{3.26}), and then derives (\ref{3.34}).
Owing to Fatou's lemma, for every $p\geq 1$, any $R>0$, it follows that
 \begin{eqnarray*}
{\mathbb E}\sup_{t\in [0,T]}\|\rho(t)\|_{BV_x(B_R)}^p \leq
{\mathbb E}\sup_{(t,x)\in [0,T]\times B_R}|\nabla_xX^{\tilde{v}}(t,0,x)|^p
e^{C_0Tp}\|\rho_0\|_{BV_x}^p\leq Ce^{C_0Tp}\|\rho_0\|_{BV_x}^p,
\end{eqnarray*}
with $C_0$ is given in (\ref{3.19}).

On the other hand, if we define $\rho^1(t,x)=\rho(t,x+B(t))$, then the unique stochastic entropy solution of (\ref{1.6}) meets
\begin{eqnarray}\label{5.1}
\left\{
  \begin{array}{ll}
\partial_t\rho^1(t,x)+b_1(x)\cdot\nabla_x\int_{\mathbb{R}}
f^\prime(v)u^1(t,x,v)dv
=0,   \\  u^1(t,x,v)=\chi_{\rho^1(t,x)}(v), \\  \rho^1(t,x)|_{t=0}=\rho_0(x), \\ \ b_1(x)=b(x+B(t)).
\end{array}
\right.
\end{eqnarray}

The calculations from (\ref{3.30}) to (\ref{3.32}) use again, for every $t_1,t_2\in (0,T)$ which is sufficiently near (i.e. there exists a positive real number $\lambda>0$, such that $|t_1-t_2|<\lambda$), every $R>0$, we get
\begin{eqnarray}\label{5.2}
&&\|\rho(t_1)-\rho(t_2)\|_{L^1_x(B_R)}
\cr\cr&\leq& C(\sqrt{|t_1-t_2|\log\log(1/|t_1-t_2|)} +|t_1-t_2|)\|\rho_0\|_{BV_x}\cr\cr&\leq& C|t_1-t_2|^\alpha\|\rho_0\|_{BV_x}, \quad  {\mathbb P}-a.s.,
\end{eqnarray}
for every $\alpha\in (0,1/2)$.

For $t_1,t_2\in (0,T)$, if $|t_1-t_2|\geq \lambda$, for every $\alpha\in (0,1/2)$, then
\begin{eqnarray}\label{5.3}
\|\rho(t_1)-\rho(t_2)\|_{L^1_x(B_R)}
\leq C\|\rho_0\|_{L^1(\mathbb{R}^d)}\leq \frac{|t_1-t_2|^\alpha}{\lambda^\alpha}C\|\rho_0\|_{L^1_x}, \quad  {\mathbb P}-a.s..
\end{eqnarray}

From (\ref{5.2}) and (\ref{5.3}), we conclude that for almost all $\omega\in\Omega$, $\rho\in {\mathcal C}^\alpha([0,T];L^1(\mathbb{R}^d))$ with every $\alpha\in (0,1/2)$. $\Box$

\vskip2mm\noindent
\textbf{Remark 5.1.} Even though the present result is concerned with $b=b(x)\in {\mathcal C}_b^\alpha\cap BV_{loc}(\mathbb{R}^d;\mathbb{R}^d)$, $f\in{\mathcal C}^1(\mathbb{R})$, it is appropriate for
\begin{eqnarray}\label{5.4}
\left\{
  \begin{array}{ll}
\partial_t\rho(t,x)+\mbox{div} F(\rho)
+\partial_{x_i}\rho(t,x)\circ\dot{B}_i(t)=0, \ \ (\omega,t,x)\in \Omega\times (0,T)\times \mathbb{R}^d, \\
\rho(t,x)|_{t=0}=\rho_0(x), \ x\in \mathbb{R}^d,
  \end{array}
\right.
\end{eqnarray}
with $F\in {\mathcal C}^1(\mathbb{R};\mathbb{R}^d)$ and $\rho_0\in L^\infty\cap L^1(\mathbb{R}^d)$. And now the regularity for $x$ is global on $\mathbb{R}^d$. Precisely speaking, we have

\begin{corollary} There is a unique stochastic entropy solution $\rho$ of (\ref{5.4}). Moreover, if $\rho_0\in
BV(\mathbb{R}^d)$, then $\rho\in
L^\infty_{\omega}(L^\infty_{t}(BV_x))$ and for almost all $\omega\in\Omega$, $\rho\in {\mathcal C}^\alpha([0,T];L^1(\mathbb{R}^d))$ for every $\alpha\in
(0,1/2)$.
\end{corollary}

\vskip2mm\noindent
\textbf{Proof.} The proof is analogue of the proof of Theorem \ref{the2.2}, the main difference is that we should substitute $\mathbb{R}^d$ for $B_R$ in Step 2 in the proof of Theorem \ref{the2.1}, now we give the details.

Let $\{e_j, j=1,2,_{\cdot},d\}$ be the standard orthogonal basis of $\mathbb{R}^d$, that $u_{\varepsilon_1}$ be the unique solution of (\ref{3.1}), with initial value $\chi_{\rho_0}$ and let $u_{\varepsilon_1}^2(t,x,v)=u_{\varepsilon_1}(t,x+he_j,v)$, then
\begin{eqnarray}\label{5.5}
\left\{
  \begin{array}{ll}
\partial_tu_{\varepsilon_1}^2(t,x,v)+F^\prime(v)\cdot\nabla_x u^2_{\varepsilon_1}
+\partial_{x_i}u_{\varepsilon_1}^2\circ \dot{B}_i(t)= \frac{1}{\varepsilon_1}
\Big[ \chi_{\rho^2_{\varepsilon_1}} -u_{\varepsilon_1}^2\Big], \\
u_{\varepsilon_1}^2(t,x,v)|_{t=0}=\chi_{\rho_0(x+he_j)}(v),
\\  \rho^2_\varepsilon=\rho^2_\varepsilon(t,x)=\int_{\mathbb{R}} u_{\varepsilon_1}^2(t,x,v)dv.
\end{array}
\right.
\end{eqnarray}
By using (\ref{3.11}) and (\ref{5.5}), for almost all $\omega\in\Omega$,
\begin{eqnarray}\label{5.6}
\|u_{\varepsilon_1}-u^2_{\varepsilon_1}\|_{L^1_{x,v}}
&\leq&
\int\limits_{\mathbb{R}^{d+1}_{x,v}}  \frac{1}{\varepsilon_1}\int\limits^t_0e^{\frac{s-t}
{\varepsilon_1}}|
\chi_{\rho_{\varepsilon_1}(s,X_{t,s}^v(x))}(v)-
\chi_{\rho^2_{\varepsilon_1}(s,X^v_{t,s}(x))}(v)|dsdxdv \cr\cr &&+
\int\limits_{\mathbb{R}^{d+1}_{x,v}}e^{-\frac{t}{\varepsilon_1}}
|\chi_{\rho_0(X^v_{t,0}(x))}(v)
-\chi_{\rho_0(X^v_{t,0}(x)+he_j)}(v)|dxdv\cr\cr
&\leq&
\int\limits_{\mathbb{R}^{d+1}_{x,v}}  \frac{1}{\varepsilon_1}\int\limits^t_0e^{\frac{s-t}
{\varepsilon_1}}|
\chi_{\rho_{\varepsilon_1}(s,x)}(v)-
\chi_{\rho_{\varepsilon_1}(s,x+he_j)}(v)|dsdxdv \cr\cr &&+
\int\limits_{\mathbb{R}^{d+1}_{x,v}}e^{-\frac{t}{\varepsilon_1}}|\chi_{\rho_0(x)}(v)-
\chi_{\rho_0(x+he_j)}(v)|dxdv\cr\cr &\leq&
\frac{1}{\varepsilon_1}\int\limits^t_0e^{\frac{s-t}
{\varepsilon_1}}\|u_{\varepsilon_1}(s)-u^2_{\varepsilon_1}(s)\|_{L^1_{x,v}}ds
 +e^{-\frac{t}{\varepsilon_1}}\|\rho_0(x)-\rho_0(x+he_j)\|_{L^1_x}.
\end{eqnarray}

Set $H(t)=\|u_{\varepsilon_1}(t)-u^2_{\varepsilon_1}(t)\|_{L^1_{x,v}}$,
then from (\ref{5.6}), one derives
\begin{eqnarray*}
H(t)\leq (1-e^{-\frac{t}{\varepsilon_1}})\max_{0\leq s\leq t}H(s)+e^{-\frac{t}{\varepsilon_1}}H(0), \quad {\mathbb P}-a.s.,
\end{eqnarray*}
which implies
\begin{eqnarray*}
\|u_{\varepsilon_1}(t)-u^2_{\varepsilon_1}(t)\|_{L^1_{x,v}}
\leq
\|\rho_0(x)-\rho_0(x+he_j)\|_{L^1_x}, \quad {\mathbb P}-a.s..
\end{eqnarray*}
Therefore, we achieve that
\begin{eqnarray}\label{5.7}
\|\rho_{\varepsilon_1}(t)\|_{BV_x}\leq \|\rho_0\|_{BV_x}.
\end{eqnarray}

Similarly, one gets analogues inequalities of (\ref{5.2}) and (\ref{5.3}),
\begin{eqnarray}\label{5.8}
\|\rho(t_1)-\rho(t_2)\|_{L^1_x}\leq C|t_1-t_2|^\alpha\|\rho_0\|_{BV_x}, \quad  {\mathbb P}-a.s.,
\end{eqnarray}
for every $\alpha\in (0,1/2)$. From (\ref{5.7} and (\ref{5.8}), the proof is finished. $\Box$

\section{Examples on non-existence and remarks}\label{sec6}
\setcounter{equation}{0}
In this section, we first prove Theorem \ref{the2.3} and then give some concluding remarks.

\vskip2mm\noindent
\textbf{Proof of Theorem \ref{the2.3}.}  The existence and uniqueness for weak solutions can be seen in Ambrosio \cite{Amb}, we skip it. It is sufficient to show the non-existence of $BV_{loc}$ solutions. For simplicity, let us suppose $d=2$. Our construction comes from \cite{WDGL}. Let $x,y\in\mathbb{R}$, we set $b_1(x)$ and $b_2(y)$ by
\begin{eqnarray*}
b_1(x)=1_{[0,1]}(x)x^{\frac{1}{2}}+1_{(1,\infty)}(x) x^{-\frac{1}{2}},
b_2(y)=1_{[0,\infty)}(y)
\frac{y}{1+y^2}.
\end{eqnarray*}
Then $b_1,b_2\in W^{1,1}_{loc}(\mathbb{R})\cap {\mathcal C}_b^{1/2}(\mathbb{R})$ and $-1/8\leq b_2^\prime\leq 1$, but $\sup b_1^\prime=\infty$.

We define $b(x,y)=(0,b_1(x)b_2(y))$, then $b\in W^{1,1}_{loc}(\mathbb{R}^2)\cap {\mathcal C}_b^{1/2}(\mathbb{R}^2)$,
$\mbox{div} b(x,y)=b_1(x)b_2^\prime(y)$ and $\mbox{div} b\in [-1/8,1]$. Consider the ODE below
\begin{eqnarray*}
\frac{d}{dt}X(t)=0, \ \frac{d}{dt}Y(t)=b_1(X(t))b_2(Y(t)),
\ X(0)=x\geq 0, \ Y(0)=y\geq 0,
\end{eqnarray*}
we gain
\begin{eqnarray}\label{6.1}
X(t,x)=x, \ Y(t,x,y)=g^{-1}(g(y)e^{2b_1(x)t}),
\end{eqnarray}
where $g(y)=e^{y^2}y^2 \ (y\geq 0)$, $g^{-1}$ is the inverse of $g$.

From (\ref{6.1}), for every $R>0$, every $t\geq 0$,
$(X(t),Y(t))([0,R]\times [0,R])\supset [0,R]\times [0,R]$. Moreover, we have
\begin{eqnarray}\label{6.2}
\Big(\frac{\partial(X(t),Y(t))}{\partial(x,y)}\Big)^{-1}(x,y)=
\left(  \begin{array}{cc} 1 & -(g^\prime(y))^{-1}g(y)2b_1^\prime(x)t \\
 0 & (g^\prime(g(y)e^{2b_1(x)t})(g^\prime(y))^{-1}e^{-2b_1(x)t} \\
\end{array}  \right).
\end{eqnarray}

Observing that, the unique weak solution of (\ref{2.26}) is given by  $\rho(t,x,y)=\rho_0((X,Y)^{-1}(t,x,y))$,  for every $t>0$, every $R>0$,  by (\ref{6.2}) one ends up with
\begin{eqnarray}\label{6.3}
&&\int_{[-R,R]\times [-R,R]}|\nabla_{x,y}(\rho_0((X,Y)^{-1}(t,x,y)))|dxdy\cr\cr&\geq&
\int_{[0,R]\times [0,R]}|\nabla_{x,y}\rho_0(x,y)|
\Big\|\Big(\frac{\partial(X,Y)}{\partial(x,y)}\Big)^{-1}\Big
\|\exp(\int^t_0\mbox{div}
b(X(r,x),Y(r,x,y))dr)dxdy\cr\cr&\geq&
\exp(-\frac{t}{8})
\int_{[0,R]\times [0,R]}|\nabla_{x,y}\rho_0(x,y)||(g^\prime(y))^{-1}g(y)2b_1^\prime(x)t|^3
dxdy\cr\cr&\geq&\exp(-\frac{t}{8})t
\int_{[0,R]\times [0,R]}|\nabla_{x,y}\rho_0(x,y)|\Big(\frac{y}{1+y}\Big)|b_1^\prime(x)|dxdy,
\end{eqnarray}
where in the last inequality we have used
$$
(g^\prime(y))^{-1}g(y)=\frac{y}{2(1+y)}, \quad \forall \ y\geq 0.
$$
If one chooses $\rho_0(x,y)=\rho_{0,1}(x)\rho_{0,2}(y)$ such that $\rho_{0,1},\rho_{0,2}\in BV(\mathbb{R})$ and
$\rho_{0,1}^\prime(x)\approx x^{-1/2}$ near $0+$, from (\ref{6.3}), then for every $t>0$,
\begin{eqnarray*}
\int_{[-R,R]\times [-R,R]}|\nabla_{x,y}(u_0((X,Y)^{-1}(t,x,y)))|dxdy
\geq C\exp(-\frac{t}{8})t
\int_0^{\epsilon}x^{-\frac{1}{2}}|b_1^\prime(x)|^3dx =\infty,
\end{eqnarray*}
where $\epsilon>0$ is a small enough positive real number. From this we complete the proof. $\Box$

\vskip2mm\noindent
\textbf{Remark 6.1.} Our existence, uniqueness and regularity results can be extended to the case of vector field $b$ is time dependent if one assumes the boundedness of $b$ on time in addition. Now Theorem \ref{the2.1} and Theorem \ref{the2.2} hold ad hoc for  $f(\rho)=\rho$ :
\begin{eqnarray}\label{6.4}
\left\{
\begin{array}{ll}
\partial_t\rho(t,x)+b(t,x)\cdot\nabla\rho(t,x)+\partial_{x_i}\rho(t,x)\circ \dot{B}_i(t)=0, \ \ (\omega,t,x)\in \Omega\times(0,T)\times \mathbb{R}^d, \\
\rho(t,x)|_{t=0}=\rho_0(x), \ x\in \mathbb{R}^d.
 \end{array}
\right.
\end{eqnarray}
From Theorem 25 \cite{FGP}, to make Theorem \ref{the2.1} and Theorem \ref{the2.2} hold for (\ref{6.4}), H\"{o}lder regularity of $b$ is enough. So an interesting question is raised: without $BV_{lov}$ regularity on $b$, do the Theorem \ref{the2.1} and Theorem \ref{the2.2} hold as well ?
Our present results do not cover this case.

\vskip2mm\noindent
\textbf{Remark 6.2.} When the stochastic perturbation vanishes in (\ref{6.4}), it becomes:
\begin{eqnarray}\label{6.5}
\left\{
  \begin{array}{ll}
\partial_t\rho(t,x)+b(t,x)\cdot\nabla\rho(t,x)=0, \ \ (t,x)\in (0,T)\times \mathbb{R}^d, \\
\rho(t,x)|_{t=0}=\rho_0(x), \ x\in \mathbb{R}^d.
  \end{array}
\right.
\end{eqnarray}
Another interesting question posed by Crippa and De Lellis \cite{CD} (similar question can be seen in \cite{Bre}) is that :

\textbf{Question.} Is there a complete topological vector space $S(\mathbb{R}^d)$ such that
$$
BV_{loc}\cap L^\infty(\mathbb{R}_t\times\mathbb{R}^d_x)\subset S(\mathbb{R}_t\times\mathbb{R}^d_x)\subset L^1_{loc}\cap L^\infty(\mathbb{R}_t\times\mathbb{R}^d_x)
$$
with the following properties ?

$\bullet$ The topology of $S(\mathbb{R}\times\mathbb{R}^d)$ is finer than the topology of $L^1_{loc}\cap L^\infty$ and coarser than the topology of $BV_{loc}\cap L^\infty$;

$\bullet$ Bounded subsets of $S(\mathbb{R}\times\mathbb{R}^d)$ are relatively compact in $L^1_{loc}\cap L^\infty$, here $\rho_k\rightarrow \rho$ in $L^1_{loc}\cap L^\infty$ if $\|\rho_k\|_{L^\infty(K)}$ is uniformly bounded and $\|\rho_k-\rho\|_{L^1(K)}$
converges to 0 for every open set $K\subset\subset \mathbb{R}^d$;

$\bullet$ If $\rho_0\in S(\mathbb{R}^d)$, $b\in S^d(\mathbb{R}\times\mathbb{R}^d)$, and $\mbox{div}_x b\in L^\infty(\mathbb{R}_t\times\mathbb{R}^d_x)$, then there exists
a (possibly non unique) $\rho\in L^\infty(\mathbb{R}^+; S(\mathbb{R}^d))$ which solves (\ref{6.5}) is the sense of distributions.

By using Depauw's construction (see \cite{Dep}), in \cite{CD}, the authors gave a negative answer for above question to $d\geq 3$. Under the stochastic perturbation of Brownian type, (\ref{6.4}) and (\ref{1.6}) may be well-posed for some particular space $S(\mathbb{R}_t\times\mathbb{R}^d_x)$, but new ideas and methods are needed to approach this problem.

\vskip6mm\noindent
\textbf{\large{Acknowledgements}}
 \vskip4mm\par
This research was partly supported by the NSF of China grants 11501577 and 11301146.

\end{document}